\documentclass{conm-p-l}

\newtheorem{theorem}{Theorem}[section]
\newtheorem{lemma}[theorem]{Lemma}
\newtheorem{prop}[theorem]{Proposition}
\newtheorem{cor}[theorem]{Corollary}

\theoremstyle{definition}
\newtheorem{definition}[theorem]{Definition}

\theoremstyle{remark}
\newtheorem{remark}[theorem]{Remark}

\numberwithin{equation}{section}

\usepackage{amsmath}
\usepackage{amsfonts}
\usepackage{amssymb}

\newcommand{\abs}[1]{\ensuremath{\left| #1 \right| }}
\newcommand{\qPs}{$q$--Pochhammer symbol}

\newcommand{\bhs}{basic hypergeometric series}
\newcommand{\rhs}{right hand side}
\newcommand{\lhs}{left hand side}

\newcommand{\wrt}{with respect to}
\newcommand{\od}{one dimensional}

\newcommand{\BT}{Bailey Transform}
\newcommand{\BL}{Bailey Lemma}

\newcommand{\RRis}{Rogers--Ramanujan identities}

\newcommand{\hgr}{hyperoctahedral group}

\newcommand{\Js}{Jackson sum}

\newcommand{\ci}{cocycle identity}

\newcommand{\tns}{\ensuremath{_{10}\varphi_9}}

\newcommand{\W}[1]{\ensuremath{\mathbb{W}(#1)}}
\newcommand{\oW}[1]{\ensuremath{\mathfrak{W}(#1)}}

\newcommand{\Jc}{Jackson coefficients}
\newcommand{\eJc}{elliptic Jackson coefficients}
\newcommand{\Mf}{Macdonald functions}
\newcommand{\eMf}{elliptic Macdonald functions}

\begin{document}

\title{Well--Poised Macdonald Functions $W_\lambda$ and \\
Jackson Coefficients $\omega_\lambda$ On $BC_n$}  

\author{Hasan Coskun}
\address{Department of Mathematics, Texas A\&M
  University--Commerce, Binnion Hall, Room 314, Commerce, TX 75429}  
\curraddr{Department of Mathematics, Texas A\&M
  University--Commerce, Binnion Hall, Room 314, Commerce, TX 75429}  
\email{hasan\_coskun@tamu-commerce.edu}

\author{Robert A. Gustafson}
\address{Department of Mathematics, Texas A\&M
  University, Milner Hall, Room 208, College Station, TX 77843}
\curraddr{Department of Mathematics, Texas A\&M
  University, Milner Hall, Room 208, College Station, TX 77843}
\email{rgustaf@math.tamu.edu}

\date{November 30, 2004} % and, in revised form, June 22, 1994.}


\begin{abstract}
The very well--poised \eMf\ $W_{\lambda/\mu}$ %(z;q, p, t, a,b)$ 
in $n$ independent variables %$z=(x_1,\ldots, x_n) \in \mathbb{C}^n$ 
are defined and their properties are
investigated. The $W_{\lambda/\mu}$ are generalized by introducing an
extra parameter %$r\in\mathbb{C}$ 
to the \eJc\ $\omega_{\lambda/\mu}$ % (z;r,q, p, t; a,b)$ 
and their properties are studied. $BC_n$ multivariable
Jackson sums in terms of both $W_{\lambda}$ and $\omega_\lambda$
functions are proved. 
\end{abstract}

\maketitle

\section{Introduction}
\label{section1}
A family of
symmetric elliptic functions $W_\lambda$ associated to the root system
$BC_n$ are introduced which extend the
symmetric Macdonald polynomials~\cite{Macdonald5}  and
Okounkov's~\cite{Okounkov1}  symmetric interpolation
polynomials. These $W_\lambda$ functions satisfy a number of
identities which generalize important identities for classical
very-well-poised basic hypergeometric series, such as the Jackson or
$q$-Dougall summation theorem and Bailey's $_{10}\phi_9$
transformation. The symmetric elliptic
$\omega_\lambda$ function which extends $W_\lambda$ by an additional
parmeter is also defined in this paper. 
However, the basic (trigonometric) version of $\omega_\lambda$
functions, together with the basic \tns\ transformation are first
given in~\cite{Coskun}. Some of these results 
are proved independently and by different methods by
Rains~\cite{Rains1}.   

The genesis of the  $W_\lambda$ functions goes back 
to a family of  symmetric  interpolation polynomials defined by
Biedenharn and Louck~\cite{BiedenharnL1}.  The proof given by Biedenharn
and Louck for the symmetry of their polynomials depends on a
well-known hypergeometric transformation formula for balanced
terminating $_4F_3$ hypergeometric series. By using Bailey's
$_{10}\phi_9$ transform, it was possible to include two additional
parameters to their functions in the two variable case, obtaining
the two variable (trigonometric) $W_\lambda$ functions.  The general
$W_\lambda$  functions were then studied as a means to generalize these
classical hypergeometric methods to the context of multivariate
hypergeometric series.  By using the elliptic version of the Bailey's
$_{10}\phi_9$ transformation as given in Frenkel and
Turaev~\cite{FrenkelT1}, it was easily seen that the same approach
generalized to the elliptic case. This extension of the trigonometric
$W_\lambda$ and $\omega_\lambda$ to the elliptic case is 
conjectured first in~\cite{Coskun}.  Following Frenkel and Turaev,
there has been much work done recently in studying different
generalizations of hypergeometric series involving elliptic functions, e.g.
\cite{Rosengren1}, \cite{SpiridonovZ1} and \cite{Warnaar2}. 

The initial motivation for the basic
(trigonometric) version of $\omega_{\lambda}$ was to prove the
$BC_n$ \BL s given in~\cite{Coskun}. \BL\ is a powerful iterative
method classically used for proving \bhs\ identities (see
\cite{AgarwalA1}, \cite{Bailey1}, \cite{AndAskRoy}, \cite{GarretI1},
\cite{LillyM3},  \cite{Spiridonov1}, \cite{Warnaar1}, \cite{Coskun1}
and \cite{Coskun2} for example). However, as 
described in Theorem~\ref{thm:bigInduction}, properties of each
function $W_\lambda$ or $\omega_{\lambda}$ are used in the proof of
the properties of the 
other. For example, that $W_\lambda(x_1, \ldots, x_n; q, p, t, a, b)$
is symmetric follows from the symmetry of $\omega_\lambda(x_1, \ldots,
x_n; r, q, p, t; a, b)$ in $x_i$ variables.

The paper is organized as follows.
At the start of \S\ref{section2} we give a definition (Definition
\ref{interpdef})  characterizing the 
symmetric elliptic $W_\lambda$ functions as solutions of a
multivariate interpolation problem.  This definition is a
analogous to that given by Okounkov~\cite{Okounkov2} for
the symmetric interpolation polynomials associated to the root system
$C_n$. In \S\ref{combindef} , a combinatorial definition of
the elliptic $W_\lambda$ functions is given  in  definitions
(\ref{definitionHfactor}), (\ref{definitionSkewW}), and
(\ref{eqWrecurrence}). Some of the basic properties of these
combinatorial $W_\lambda$ functions is also proved.  In
\S\ref{multiseries}, a definition (Def. \ref{multipleBHS}) of multiple
very-well-poised elliptic 
hypergeometric series is given.  The definition of the elliptic
$\omega_\lambda$ or more generally the $\omega_{\lambda/\nu}$
functions~(\ref{eq:omega{lambda,mu}}) is also given in this section
and some of their 
properites are proved. In particular, $\omega_{\lambda}$ fuction is
characterized in terms of certain multiplication and shift operators,
analogous to the characterization of Biedenharn and Louck's factorial
Schur functions given by Goulden and Hamel~\cite{GouldenH1}. It is
also shown that $\omega_{\lambda}$ is independent of different
representations of $\lambda$.  

In \S\ref{section3} Theorem~\ref{thm:bigInduction} is proved by
induction in nine steps. The first steps (\ref{def:item1})--
(\ref{def:item4}) establish that the combinatorial $W_\lambda$
functions satisfy the conditions given in Definition
\ref{interpdef}.  Other steps include the proof of the $W$--Jackson sum
(\ref{eq:W_Jackson}), a $BC_n$ generalization of elliptic Bailey's \tns\
transformation (\ref{BCn10phi9t}), and an elliptic cocycle identity
for the $\omega$ functions.  There is also a discussion of a $BC_n$
generalization of the Bailey Transform, which is a key algorithm used
to generate classical hypergeometric series identities
(see~\cite{Andrews}, \cite{AndAskRoy}, \cite{Bailey2}, \cite{LillyM3},
\cite{Coskun1} and \cite{Coskun2}). Last few
results mentioned here are proved by easily extending the methods
developed for the basic (trigonometric) case 
in~\cite{Coskun} to the elliptic case.

Finally, \S\ref{section4} gives a proof of the
$\omega$--\Js~(\ref{10phi9sum}) which extends the $W$--Jackson sum by
an additional  
parmeter $r\in\mathbb{C}$. The proof uses the symmetry and vanishing
properties of $\omega_\lambda$, and a key
identity~(\ref{omega(xr^delta)}) for  
 $\omega_{\lambda/\mu}(xr^{\delta(m)};r, q, p, t; a,b) $.
The $\omega$--\Js\ reduces to the $W$--\Js\ in the limit as
$r\rightarrow t$. The basic (trigonometric) case of
$\omega$--\Js~(\ref{10phi9sum}) is proved 
in~\cite{Coskun}. We finish this section by listing   
further important properties of the  $\omega_\lambda$ function. 

\section{Elliptic Macdonald Functions And Elliptic \Jc}
\label{section2}
We recall some terminology from the theory of basic hypergeometric
series. For $a,q\in \mathbb{C}$, $|q|<1$ and integer $m$, we define
the basic factorial $(a;q)_\infty:=\prod_{i=0}^{\infty} (1-aq^i)$ and
$(a;q)_m :=(a;q)_\infty/(aq^m;q)_\infty$. We also define an elliptic
analogue of the basic factorial as follows. For $x, p\in \mathbb{C}$
and $\abs{p}<1$, let the elliptic function be given by
\begin{equation}
  E(x) = E(x;p) := (x; p)_\infty (p/x; p)_\infty
\end{equation}
and for $a\in \mathbb{C}$, and a positive integer $m$ define
\begin{equation}
  (a; q,p)_m := \prod_{k=0}^{m-1} E(aq^m)
\end{equation}
The definition is extended to negative $m$ by setting
  $(a; q,p)_m = 1/ (aq^{m}; q, p)_{-m} $. 
Note also that when $p=0$, $(a; q,p)_m$ reduces to standard
(trigonometric) \qPs.

For any partition $\lambda = (\lambda_1, \ldots, \lambda_n)$ and
$t\in\mathbb{C}$ we can also define
\begin{equation}
\label{ellipticQtPocSymbol}
  (a)_\lambda=(a; q, p, t)_\lambda := \prod_{k=1}^{n} (at^{1-i};
  q,p)_{\lambda_i} .
\end{equation}
Note that when $\lambda=(\lambda_1) = \lambda_1$ is a single part
partition, then $(a; q, p, t)_\lambda = (a; q, p)_{\lambda_1} =
(a)_{\lambda_1}$. We'll also use the notation
\begin{equation}
  (a_1, \ldots, a_k)_\lambda = (a_1, \ldots, a_k; q, p, t)_\lambda :=
  (a_1)_\lambda \ldots (a_k)_\lambda .
\end{equation}

Let $x_1, \ldots, x_n\in \mathbb{C}$ be a set of variables and
$\mathbb{F} = \mathbb{C}(q,p,t,a,b)$ field of rational functions in
the parameters $q,p,t,a,b\in \mathbb{C}$. For positive integers $i,k$
with $1\leq i\leq n$, let $V(i, k)$ be the vector space of functions
in $x_1, \ldots, x_n$ spanned over $F$ by $1$ and
\[\{1/E(bq^m x_i), 1/E(bq^ma^{-1}x_i^{-1}) | 1\leq m \leq k \}\]
Set $V(i, 0) = \mathbb{F}$ and $V(i) = \cup_{k\geq 0} V(i,k)$. Also
let $V_n[k]$ be the vector space over $\mathbb{F}$ generated by the
product $V(1,k)\,V(2,k) \cdots V(n,k)$ and $V_n$ the vector space
generated by $V(1)\,V(2) \cdots V(n)$.

Let $\mathcal{S}$ be the hyperoctahedral group of symmetries generated
by the permutations of the variables $x_1, \ldots, x_n$ together with
the inversions $x_i \rightarrow 1/ax_i$ for $1\leq i\leq n$. Then
$V^\mathcal{S}_n[k]$ (respectively $V^\mathcal{S}_n$) is the subspace
of $V_n[k]$ $(V_n)$ invariant under $\mathcal{S}$.

At various times it may be necessary to assume that $p,q,t,a,b$ do not
take on special values or that $\abs{t}<1$. It will generally be clear
when this is required.

Note the effect on $V(i, k)$ of setting $p=0$. Then $V(i,k)$ reduces
to a subspace of the field of rational functions
$\mathbb{C}(x_1,\ldots, x_n)$ over $\mathbb{C}(q,t,a,b)$ spanned by
$1$ and
\[\{ 1/(1-bq^m x_i), 1/(1-bq^ma^{-1}x_i^{-1}) | 1\leq m \leq k \}\]
After setting $p=0$, we may also consider the limits $b\rightarrow
\infty$ or $b\rightarrow 0$. The effect on $V(i)$ is essentially to
shift the poles of $1/(1-bq^m x_i)$ and $1/(1-bq^ma^{-1}x_i^{-1})$ to
0 or $\infty$. Consequently, the elements of $V^\mathcal{S}_n$ will
asymptotically tend toward symmetric polynomials in the variables
$x_i^{-1}$ and $ax_i$.

We will now characterize a special basis of $V^\mathcal{S}_n$. This
basis will extend, up to a change of variables and normalization, the
$BC_n$ type interpolation polynomials $P_\lambda^{*}(x;q,t,s)$ of
Okounkov~\cite{Okounkov1} and consequently also extend the shifted
Schur functions of F. Knop~\cite{Knop1}, A. Okounkov~\cite{Okounkov2},
G. Olshanski~\cite{OkOl} and S. Sahi~\cite{Sahi1}. These in turn
extend the homogeneous Macdonald polynomials
$P_\lambda(x;q,t)$~\cite{Macdonald1} and Schur functions.

\begin{definition}
\label{interpdef}
Let $t^{\delta(n)}=(t^{n-1}, t^{n-2}, \ldots, 1)$ and for $n$-part
partition $\lambda$ $q^\lambda=(q^{\lambda_1}, \ldots,
q^{\lambda_n})$. Also define $-\lambda^r = (-\lambda_n,
-\lambda_{n-1}, \ldots, -\lambda_1)$ and set $q^{-\lambda^r}t^{\delta(n)}
= (q^{-\lambda_n}t^{n-1}, q^{-\lambda_{n-1}}t^{n-2}, \ldots,
q^{-\lambda_1})$. Define $W_\lambda(x_1, \ldots, x_n; p, q, t, a, b)$
to be the element of $V^\mathcal{S}_n$ satisfying the following
conditions:
%\begin{enumerate}
%\item
~\begin{equation}
\label{def:item1}
W_\lambda(x; q,p, t, a, b) \in V^\mathcal{S}_n[\lambda_1]
\end{equation}
%\item
\begin{equation}
\label{def:item2}
W_\lambda(q^\nu t^{\delta(n)}; q,p, t, a, b) =0
\end{equation}
if $\lambda\not\subseteq \nu$.
\item Setting
\begin{equation}
\tilde{W}_\lambda(x; q, p, t, a, b) = W_\lambda(x;
  q,p,t, a,b) \prod_{i=1}^n \dfrac{(qbx_i, qba^{-1}x_i^{-1})_{\lambda_1}
  } {(x_i^{-1}, ax_i)_{\lambda_1} } .
\end{equation}
Then
\begin{equation}
\label{def:item3}
\tilde{W}_\lambda(b^{-1}
  q^{-\nu^r} t^{\delta(n)}; q, p, t, a, b) =0
\end{equation}
if $\nu\not\subseteq \lambda$ and $\nu\subseteq (\lambda_1)^n$.
%\item
\begin{equation}
\label{def:item4}
W_\lambda(q^\lambda t^{\delta(n)}; q,p, t, a, b) = N(\lambda,
  n; q,p, t, a, b)
\end{equation}
is a normalization constant to be specified below.
%\end{enumerate}
\end{definition}

Neither the existence nor the uniqueness of the $W_\lambda$ are
obvious. We will construct the $W_\lambda$ functions below satisfying
conditions (1)--(4). The set of the $W_\lambda$ functions as $\lambda$
ranges over all $n$-part partitions $\lambda$ will form a basis of
$V^\mathcal{S}_n$ over $\mathbb{F}$. The uniqueness of the $W_\lambda$
functions will then follow from conditions (2)--(4).

\begin{remark}
The definition of the $W_\lambda$ functions is similar to to
Okounkov's $P^\star_\lambda$ functions~\cite{Okounkov1}. The
difference is that conditions (1) and (3) above are replaced by a
multidegree condition $\deg P_\lambda^\star(x; q,t,s) \leq |\lambda|$
or $\deg P_\lambda^\star(x; q,t) \leq |\lambda|$ where $\abs{\lambda} =
\lambda_1+\ldots +\lambda_n$ is the weight of the $n$-part partition.
\end{remark}

The normalization we will use for the $W_\lambda$ functions is given
as follows. Let $\lambda=(\lambda_1, \ldots, \lambda_n)$ be a
partition, then
\begin{multline}
N(\lambda, n; q,p, t, a, b) = W_\lambda(q^\lambda t^{\delta(n)}; q,p,
t, a, b) \\
= \prod_{k=1}^n\left\{\frac{(qbt^{n-k}, qt^{n-k})_{\lambda_k}
    (at^{2n-2k})_{2\lambda_k}}
{((a/b)t^{n-k}, at^{n-k})_{\lambda_k} (qbt^{n+1-2k})_{2\lambda_k}}
t^{(n+1-2k)\lambda_k}\right\}\\
\cdot(a/(qb))^{|\lambda|}\cdot\prod_{1\leq i < j\leq
n}\frac{(qt^{j-i-1})_{\lambda_i-\lambda_j}(at^{2n-i-j})_{\lambda_i+\lambda_j}}
{(qt^{j-i})_{\lambda_i-\lambda_j}(at^{1+2n-i-j})_{\lambda_i+\lambda_j}}
\end{multline}

\subsection{Combinatorial Definition of $W_\lambda$ Functions.}
\label{combindef}
For positive integer $n$, let $\lambda=(\lambda_1, \ldots,
\lambda_n)$,  $\mu=(\mu_1, \ldots, \mu_n)$ be partitions such that the
skew partition $\lambda/\mu$ is a horizontal strip; i.e. $\lambda_1
\geq \mu_1 \geq\lambda_2 \geq \mu_2 \geq \ldots \lambda_n \geq
\mu_n$. Also set $\lambda_{n+1} = \mu_{n+1} =0$. For $\lambda$ and
$\mu$ as above and $b\in\mathbb{C}$, we define
\begin{multline}\label{definitionHfactor}
H_{\lambda/\mu}(q,p,t,b) \\
:= \prod_{1\leq i < j\leq
n}\left\{\dfrac{(q^{\mu_i-\mu_{j-1}}t^{j-i})_{\mu_{j-1}-\lambda_j}
(q^{\lambda_i+\lambda_j}t^{3-j-i}b)_{\mu_{j-1}-\lambda_j}}
{(q^{\mu_i-\mu_{j-1}+1}t^{j-i-1})_{\mu_{j-1}-\lambda_j}(q^{\lambda_i
    +\lambda_j+1}t^{2-j-i}b)_{\mu_{j-1}-\lambda_j}}\right.\\
\left.\cdot 
\dfrac{(q^{\lambda_i-\mu_{j-1}+1}t^{j-i-1})_{\mu_{j-1}-\lambda_j}}
{(q^{\lambda_i-\mu_{j-1}}t^{j-i})_{\mu_{j-1}-\lambda_j}}\right\}\cdot\prod_{1\leq
i <(j-1)\leq n}
\dfrac{(q^{\mu_i+\lambda_j+1}t^{1-j-i}b)_{\mu_{j-1}-\lambda_j}}
{(q^{\mu_i+\lambda_j}t^{2-j-i}b)_{\mu_{j-1}-\lambda_j}}.
\end{multline}
and also for $x\in \mathbb{C}$,
\begin{multline}\label{definitionSkewW}
W_{\lambda/\mu}(x; q,p,t,a,b)
:= H_{\lambda/\mu}(q,p,t,b)\cdot\dfrac{(x^{-1}, ax)_\lambda
  (qbx/t, qb/(axt))_\mu}
{(x^{-1}, ax)_\mu (qbx, qb/(ax))_\lambda}\\
\cdot\prod_{i=1}^n\left\{\dfrac{E(bt^{1-2i}q^{2\mu_i})}{E(bt^{1-2i})}
  \dfrac{(bt^{1-2i})_{\mu_i+\lambda_{i+1}}}
{(bqt^{-2i})_{\mu_i+\lambda_{i+1}}}\cdot
t^{i(\mu_i-\lambda_{i+1})}\right\},
\end{multline}

Note that for $t=q$, $H_{\lambda/\mu}(q,p,q,b)=1$, and for $p=b=0$,
$H_{\lambda/\mu}(q,0,t,0)=\psi_{\lambda/\mu}$ where
$\psi_{\lambda/\mu}$ is the weight function in Macdonald's
combinatorial formula~\cite{Macdonald1} :
\begin{equation}
P_\lambda(x_1, \ldots, x_n) = \sum_{\mu \prec\lambda} \psi_{\lambda/\mu}
x_1^{\abs{\lambda/\mu } } P_\mu(x_2, \ldots, x_n).
\end{equation}
where $\mu \prec\lambda$ means
$\lambda_1 \geq \mu_1 \geq\lambda_2 \geq \mu_2 \geq \ldots \lambda_n \geq
\mu_n\geq 0$. (In this case, we can also assume $\mu_n=0$, since
$P_{\mu}(x_2, \ldots, x_n) =0$ otherwise.)

Letting $\lambda$, $\mu$ be arbitrary partitions $\lambda=(\lambda_1, 
\ldots,
\lambda_n)$, $\mu=(\mu_1, \ldots, \mu_n)$ with $\lambda/\mu$ a skew
partition, we define the function $W_{\lambda/\mu}(y, z_1, \ldots,
z_\ell; q,p,t,a,b)$ in $\ell+1$ variables $y, z_1, \ldots, z_\ell
\in\mathbb{C}$ by the following recursion formula

\begin{multline}
\label{eqWrecurrence}
W_{\lambda/\mu}(y,z_1,z_2,\ldots,z_\ell;q, p, t, a, b) \\
= \sum_{\nu\prec \lambda} W_{\lambda/\nu}(yt^{-\ell};q, p, t, at^{2\ell},
bt^\ell) \, W_{\nu/\mu}(z_1,\ldots, z_\ell;q, p, t, a, b).
\end{multline}

We also set $W_{\lambda/0}(x;q, p, t,a, b) = W_{\lambda}(x;q, p, t, a,
b)$ where $x\in\mathbb{C}^n$.

Note that when $\lambda$ has only one part and $x\in \mathbb{C}$, then
\begin{equation}
W_{\lambda}(x; q,p,t,a,b)
:=\dfrac{(x^{-1}, ax)_\lambda }
{(qbx, qb/(ax))_\lambda}
\end{equation}
is independent of the parameter $t$.

Notice also the strange property that for any $n$-part partition
$\lambda\neq 0$, and $x\in\mathbb{C}$, $W_{\lambda/\lambda}(x;q, p, t, a,
b)$ is not identically 1. In fact, $W_{\lambda/\lambda}(x;q, p, t, a,
b)$ is not even a constant.

We'll now discuss some properties of the combinatorially defined
$W_\lambda$ functions. The first proposition is an easy consequence of
the definitions~(\ref{definitionHfactor}),~(\ref{definitionSkewW})
and~(\ref{eqWrecurrence}).

\begin{prop}
Let $\lambda$ be an $n$-part partition with $\lambda_n\neq 0$ and
$0\leq k\leq \lambda_n$ for some integer $k$. Let
$x=(x_1,\ldots,x_n)\in \mathbb{C}^n$, then
\begin{multline}
\label{reduction}
W_{\lambda}(x;q, p, t, a, b)
=\prod_{j=1}^{[n/2]} \dfrac{(qbt^{n-2j})_{2k}}
{(qbt^{-1-n+2j})_{2k}}
\prod_{i=1}^n\dfrac{(x_i^{-1})_k (ax_i)_k}{(qbx_i)_k (qb/(ax_i))_k} \\
\cdot W_{\lambda-k^n}(xq^{-k};q, p, t, aq^{2k}, bq^{2k})
\end{multline}
where $k^n$ is the $n$-part partition all of whose parts equal the
integer $k$.
\end{prop}

\begin{cor}
\begin{equation}
\label{reductionSpecial}
W_{\lambda_1^n}(x;q, p, t, a, b)
=\prod_{j=1}^{[n/2]} \dfrac{(qbt^{n-2j})_{2\lambda_1}}
{(qbt^{-1-n+2j})_{2\lambda_1}}
\prod_{i=1}^n\dfrac{(x_i^{-1})_{\lambda_1} (ax_i)_{\lambda_1} }
{(qbx_i)_{\lambda_1} (qb/(ax_i))_{\lambda_1} }
\end{equation}
\end{cor}

\begin{remark}
An alternative normalization of the $W_\lambda$ functions are given
as follows:
\begin{multline}
\label{normalizedW}
W^*_{\lambda}(x;q, p, t, a, b) \\
= \prod_{1\leq i<j \leq n} \left\{\dfrac{(t^{j-i})_{\lambda_i-\lambda_j}
(q bt^{n-i-j} )_{\lambda_i+\lambda_j}} { (t^{j-i+1})_{\lambda_i-\lambda_j}
(q b t^{n-i-j+1})_{\lambda_i+\lambda_j}}
\right\} W_{\lambda}(x;q, p, t, a, b)
\end{multline}
One advantage of this normalization is that it simplifies the above
formula. We have
\begin{multline}
\label{reduction2}
W^*_{\lambda}(x;q, p, t, a, b)
= \prod_{i=1}^n\dfrac{(x_i^{-1})_k (ax_i)_k}{(qbx_i)_k (qb/(ax_i))_k}
W^*_{\lambda-k^n}(xq^{-k};q, p, t, aq^{2k}, bq^{2k})
\end{multline}
\end{remark}

Another property is the vanishing theorem for the combinatorially
defined $W_\lambda$ function.
\begin{theorem}
\label{thm:vanishingW}
Let $\pi=(\pi_1,\ldots,\pi_n)$ and $\lambda =(\lambda_1,\ldots,\lambda_n)$ 
be
$n$-part partitions such that $\lambda \nsubseteq \pi$. We then have
\begin{equation}
\label{FundamentalVanishing}
  W_{\lambda}(q^\pi t^\delta;q, p,t,a,b)=0.
\end{equation}
\end{theorem}

\begin{proof}
By induction on $n$ and using the recursive
definition~(\ref{eqWrecurrence}) of the $W$ functions
\begin{multline}
W_{\lambda}(q^\pi t^\delta;q, p, t,a,b)\\
=\sum_\nu W_{\lambda/\nu}(y_1;q, p, t,at^{2(n-l+1)},bt^{n-l+1})\cdot
W_{\nu}(y_2;q, p t,a,b),
\end{multline}
where $y_1=(q^{\pi_1}t^{l-2},q^{\pi_2}t^{l-2},\ldots,q^{\pi_{l-1}})$ and
$y_2=(q^{\pi_l}t^{n-l},q^{\pi_{l+1}}t^{n-l-1},\ldots,q^{\pi_n})$ and
$l$, $1\leq l\leq n$ is chosen so that $\lambda_l > \pi_l$ (since
$\lambda\not\subseteq \pi$).
\end{proof}

\begin{cor}
If $\lambda = (\lambda_1, \ldots, \lambda_n)$ is a partition with
$\lambda_{n-k+1}=\ldots =\lambda_n =0$, where $1\leq
  k\leq n$, and $x=(x_1,\ldots,x_{n-k},t^{\delta(k)})\in \mathbb{C}^n$
  and $t^{\delta(k)}=(t^{k-1},\ldots, 1)$, then set
  $\hat{\lambda}=(\lambda_1, \ldots, \lambda_{n-k})$ and
  $\hat{x}=(x_1,\ldots,x_{n-k})$. We have
\begin{equation}
\label{eq:ShiftedStable}
W_\lambda(x;q,p,t,a,b) = W_{\hat{\lambda}}(\hat{x}t^{-k};q,p,t,
at^{2k}, bt^{k}).
\end{equation}
\end{cor}

\begin{proof}
Use the recurrence relation
\begin{multline}
\label{eq:PfShiftedStable}
W_{\lambda}(x;q, p, t, a, b)
= \sum_{\mu\prec \lambda} W_{\lambda/\mu}(\hat{x}t^{-k};q,p,t,
at^{2k}, bt^{k}) \, W_{\mu}(t^{\delta(k)}; q, p, t, a, b).
\end{multline}
By Theorem~(\ref{thm:vanishingW}), the only non--vanishing term
in~(\ref{eq:PfShiftedStable}) is when $\mu=0^k$. The result follows
after observing that $W_{0^k}(t^{-\delta(k)}; q, p, t, a, b) =1$ and
\[W_{\hat{\lambda}/0^k}(\hat{x}t^{-k}; q,p,t, at^{2k}, bt^{k}) =
W_{\hat{\lambda}}(\hat{x}t^{-k};q,p,t, at^{2k}, bt^{k}).\]
\end{proof}

A further property of the combinatorially defined $W_\lambda$ function
corresponds to a property of Okounkov's $BC$-type interpolation
polynomials~\cite{Okounkov1}:
\begin{multline}
P^*_{\mu}(1/x_1, \ldots, 1/x_n; 1/q, 1/t, 1/s)
= s^{2|\mu|} t^{2n-2|\mu|} P^*_{\mu}(x_1, \ldots, x_n; q, t, s).
\end{multline}
We have

\begin{prop}
\begin{multline}
W_{\lambda}(x_1^{-1}, \ldots, x_n^{-1}; q^{-1}, p, t^{-1}, a^{-1},
b^{-1})\\
= a^{-2|\lambda|} b^{2|\lambda|} q^{2|\lambda|}
t^{2n(\lambda)-2(n-1)|\lambda|} W_{\lambda}(x_1, \ldots, x_n; q, p, t,
a, b)
\end{multline}
\end{prop}

\begin{proof}
One checks that for skew $W$ functions of one variable $x\in
\mathbb{C}$, we have
\begin{multline}
W_{\lambda/\mu}(x^{-1}; q^{-1}, p, t^{-1}, a^{-1}, b^{-1})\\
= (qba^{-1})^{2|\lambda|-2|\mu|} t^{2n(\lambda)-2n(\mu) - 2|\mu|}
W_{\lambda/\mu} (x; q, p, t, a, b)
\end{multline}
and the result follows by the induction using the recursive
definition~(\ref{eqWrecurrence}).
\end{proof}

The next result is an analog for the $W_\lambda$ functions of the
reversal formulas for basic hypergeometric series.
\begin{prop}
Let $\lambda\subseteq N^n$ be a partition and $x\in\mathbb{C}^n$. Let
$\bar{\lambda} = N^n -\lambda^r = (N-\lambda_n, \ldots,
N-\lambda_1)$. Then
\begin{multline}\label{eq:reversal}
W_{\bar{\lambda}}(x; q, p, t, a, b)
= \prod_{1\leq i<j\leq n} \dfrac{(qbt^{n-i-j+1}, q^{1-N}b^{-1}
  t^{n-i-j} )_{\lambda_i+\lambda_j} } {(qbt^{n-i-j}, q^{1-N}b^{-1}
  t^{n-i-j+1} )_{\lambda_i+\lambda_j} } \\
\cdot (qba^{-1})^{2|\lambda|} \prod_{i=1}^n \dfrac{(x_i^{-1}, ax_i)_N }{
  (qbx_i, qb/(ax_i))_N }
W_{\lambda}(bq^Nx; q, p, t, ab^{-2}q^{-2N}, b^{-1}q^{-2N})
\end{multline}
which using the alternate normalization can be written as
\begin{multline}
W^*_{\bar{\lambda}} (x; q, p, t, a, b)
= W^*_{N^n}(x; q, p, t, a, b) \\
\cdot (qb/2a)^{2|\lambda|}
W^*_{\lambda}(bq^Nx; q, p, t, ab^{-2}q^{-2N}, b^{-1}q^{-2N})
\end{multline}
\end{prop}

\begin{proof}
Direct computation.
\end{proof}
We finally prove the following
\begin{prop}
$W_\lambda(x; q,p,t,a,b)$, where $x = (x_1,\ldots, x_n) = q^u= (q^{u_1}, 
\ldots, q^{u_n})\in \mathbb{C}^n $, is elliptic in each of the
variables $u_1,\ldots, u_n \in \mathbb{C}^n$ and $0<q<1$
\end{prop}
\begin{proof}
It will
suffice to show that for any variable $x_i, 1\leq i \leq n$
substituting $px_i$ in place of $x_i$ leaves the   $W_\lambda(x;
q,p,t,a,b)$ function invariant.
   From the definition (\ref{eqWrecurrence} ) of the $W$ function this
   reduces to considering the skew function in one variable 
$W_{\lambda/\mu}$.
   Using the definition of $W_{\lambda/\mu}$ given in 
(\ref{definitionSkewW})
and (\ref{definitionHfactor}),  consider the result of replacing the 
variable
$x$ by $px$  in the factor
\begin{equation}
  \label{eq:px-lambda}
  \dfrac{((px)^{-1}, apx)_\lambda}{(qbpx, qb/(apx))_\lambda} =
  \dfrac{((x)^{-1}, ax)_\lambda}{(qbx,
  qb/(ax))_\lambda}
\end{equation}
and a similarly
\begin{equation}
  \label{eq:px-mu}
\dfrac{(qbpx/t,qb/(apxt))_\mu}{((px)^{-1}, apx)_\mu} =
\dfrac{(qbx/t,qb/(axt))_\mu}{((x)^{-1}, ax)_\mu}.
\end{equation}
Hence  $W_{\lambda/\mu}$ is invariant under replacing $x$ by $px$
and  $W_\lambda(x;
q,p,t,a,b)$ is invariant under replacing $x_i$ by $px_i$.
\end{proof}
\subsection{Multiple Elliptic Hypergeometric Series}
\label{multiseries}
We'll now define an analogue of terminating very--well--poised elliptic
hypergeometric series on $BC_n$.

\begin{definition}
\label{multipleBHS}
Let $a_i$ be non--zero complex parameters for $i\in[k-1]$ where $k\geq
5$ is an integer, and $\mu$ and $\nu$ be partitions such that
$\ell(\nu)\leq n$ and $\mu\subseteq \nu$. Then set
\begin{multline}
\label{MBHS}
_{k+1}\Phi^n_k\left[a_1,a_2,\ldots, a_{k-3}; a_{k-2}, \nu \|\,
  a_{k-1}, \mu \right]_{q,p,t} \\ =
\sum_{\substack{\lambda\\  \mu \subseteq \lambda \subseteq \nu}}
   K_\lambda (a_1) \dfrac{
   [a_1 t^{1-n}, a_2, \ldots, a_{k-3}]_{\lambda}}{ [qt^{n-1}, qa_1/a_2,
   \ldots, qa_1/a_{k-3}]_{\lambda} } \\
  \cdot W_\lambda(q^\nu t^{\delta(n)}; q, p, t, a_{k-2} t^{2-2n},
   a_{1}t^{1-n}) \\ \cdot W_\mu(q^\lambda t^{\delta(n)}; q,p,t, a_{1}
   t^{2-2n}, a_{k-1}t^{1-n})
\end{multline}
where
\begin{multline}
K_\lambda(a_1) := \prod_{i=1}^{n}\left\{ \dfrac{E(a_1t^{2-2i} q^{2\mu_i})}
    {E(a_1t^{2-2i})}  \left(qt^{2i-2} \right)^{\mu_i} \right\} \\
\cdot \prod_{1\leq i<j \leq n} \left\{\dfrac{
   (qt^{j-i})_{\lambda_i-\lambda_j} (a_1
   t^{3-i-j})_{\lambda_i+\lambda_j}}
{(qt^{j-i-1})_{\lambda_i-\lambda_j} (a_1 t^{2-i-j})_{\lambda_i+\lambda_j} }
\right\}
\end{multline}
and $\subseteq$ denotes the partial inclusion ordering defined by
\begin{equation}
\label{partialordering}
\mu \subseteq \lambda \;\Leftrightarrow \;\mu_i \leq \lambda_i, \quad
\forall i\geq 1.
\end{equation}

We'll sometimes suppress the $q,p,t$ dependence. We occationally drop
$n$ in the notation as well, for the series does not depend on $n$ as
we verify below. In particular when $\mu=0$, we simplify the notation for 
the
\lhs\ of~(\ref{MBHS}) and write
\begin{multline}
\label{MBHSmu=0}
{_{k+1}\Phi_k} \left[a_1,a_2,\ldots, a_{k-3}; a_{k-2}, \nu \right] \\
:= {_{k+1}\Phi^n_k} \left[a_1,a_2,\ldots, a_{k-3}; a_{k-2}, \nu \|\,
  a_{k-1}, 0 \right]_{q,p,t}
\end{multline}
\end{definition}

Note that the series~(\ref{MBHS}) and~(\ref{MBHSmu=0}) are very
well--poised, in the sense that the
denominator parameters have the form $qa_1/a_i$, where $a_i$ are the
numerator parameters, and the very--well--poised factor is located inside the
factor $K_\lambda (a_1)$. We'll
say that the series~(\ref{MBHS}) and~(\ref{MBHSmu=0}) are balanced, if the
product of the denominator parameters is $qt^{n-1}$ times that of numerator
parameters.
\begin{remark}
The definition~(\ref{multipleBHS}) is a direct extension of the
multiple basic hypergeometric series given in~\cite{Coskun} to the
elliptic case.   
\end{remark}

Just as for the trigonometric version given in~\cite{Coskun}, the
elliptic \Jc\ $\omega_\lambda(x;r,q,p,t,a,b)$ are first defined
explicitely in the case of one variable $\omega_{\lambda/\mu}$
function associated to a skew partition $\lambda/\mu$. Then the
multivariable $\omega_{\lambda}$ function is defined in terms of
the $\Omega$ algebra generated by certain shift and multiplication
operators on the $\mathbb{Z}$--space $V$ of infinite lower--triangular
matrices defined in
\begin{definition}
\label{Valgebra}
Let $\mathbb{F}=\mathbb{C}(q, p, t,r,a,b)$ and let
$\mathbb{F}_\infty=\mathbb{F}(X)$ denote the extended field of
of rational functions over $\mathbb{F}$ in
the infinite set of indeterminates
$X=\{x_1,x_2,\ldots \}$.
We use $\mathbb{F}(z)$ for the case when $z=(x_1,\ldots, x_n)$ is
arbitrary number of variables $n\in\mathbb{Z}_{>}$. Then
$V$ denotes the $\mathbb{Z}$--space of all infinite
lower--triangular
matrices indexed by partitions,
whose entries are from $\mathbb{F}_\infty$. The condition that $u\in
V$ is lower triangular \wrt\ the partial inclusion
ordering~(\ref{partialordering}) can be phrased in the form
\begin{equation}
  u_{\lambda\mu} = 0 ,\, \quad \mathrm{when}\; \mu \not \subseteq
  \lambda
\end{equation}
The addition operation is defined in the standard form as
\begin{equation}
  (u+v)_{\lambda\mu} := u_{\lambda\mu} + v_{\lambda\mu},
\end{equation}
and the multiplication operation is defined by the relation
\begin{equation}
\label{multiplication}
  (uv)_{\lambda\mu} := \sum_{\mu\subseteq\nu\subseteq\lambda}
  u_{\lambda\nu} v_{\nu\mu}
\end{equation}
for $u,v\in V$. This operation is clearly associative since the same
matrix entries enter into the double sum either way we sum it.
\end{definition}

We now give the definition of the \eJc. Recall~\cite{Macdonald5} that
if $\mu \subseteq 
\lambda$, the set theoretic difference $\lambda-\mu$ is called a skew
(diagram) partition. A skew partition is called a
horizontal strip if it satisfies the betweenness
condition $\lambda_1 \geq \mu_1 \geq\lambda_2 \geq \mu_2 \geq \ldots $. 

\begin{definition}[Elliptic Jackson Coefficients]
\label{def:omega{lambda,mu}}
Let $\lambda$ and $\mu$ be partitions of at most $n$--parts such that
$\lambda/\mu$ is a skew partition. Then the \Jc\
$\omega_{\lambda/\mu}$ are defined by
\begin{multline}
\label{eq:omega{lambda,mu}}
\omega_{\lambda/\mu}(x; r, q,p,t; a,b)
:= \dfrac{(x^{-1}, ax)_{\lambda}} {(qbx, qb/ax)_{\lambda}}
    \dfrac{(qbr^{-1}x, qb/axr)_{\mu}}{(x^{-1}, ax)_{\mu}} \\
\cdot \dfrac{(r, br^{-1}t^{1-n})_{\mu}}{(qbr^{-2}, qt^{n-1})_{\mu}} %\\
  \prod_{i=1}^{n}\left\{ \dfrac{E(br^{-1}t^{2-2i} q^{2\mu_i})}
    {E(br^{-1}t^{2-2i})}  \left(qt^{2i-2} \right)^{\mu_i} \right\} \\
\cdot \prod_{1\leq i< j \leq n} \hspace*{-5pt} \left\{ \dfrac{
    (qt^{j-i})_{\mu_i - \mu_j} } { (qt^{j-i-1})_{\mu_i - \mu_j} }
    \dfrac{ (br^{-1}t^{3-i-j})_{\mu_i + \mu_j} } {
    (br^{-1}t^{2-i-j})_{\mu_i + \mu_j} } \right\}  %\\ \cdot
W_{\mu} (q^{\lambda}t^{\delta(n)}; q, p, t, bt^{2-2n}, br^{-1}t^{1-n})
\end{multline}
where $x,r,q,p,t, a,b\in\mathbb{C}$ are free complex parameters.
\end{definition}
Note that $\omega_{\lambda/\mu}(x; r; a,b)=\omega_{\lambda/\mu}(x; r,
q,p,t; a,b)$ is defined when the skew
partition $\lambda/\mu$ is not a horizantal strip.

Defining \W{x_i}$\,=\,$\W{x_i;r,q,p,t,a,b} for each $i\in[n]$ as the
  infinite matrix whose $(\lambda\mu)$--entry is
  $\omega_{\lambda/\mu}(x_i;r,q,p,t,a,b)$ we introduce the $\Omega$
algebra that produces the recursion formula for the \Jc\
  $\omega_\lambda$. The recurrence formula we obtain is
  completely analogous to that for $W_\lambda$ functions.

\begin{definition}
With notation as above,
let $\mathfrak{X_i}$ denote the linear operator on $V$ defined as the
right  matrix  multiplication by \W{x_i}, for each indeterminate
  $x_i$, and let $\mathfrak{I}$ denote the identity
  operator on the same space. Also define the linear operator
  $\mathfrak{M}$ on $V$ by
\begin{equation}
    (\mathfrak{M} v(X;r,q,p,t,a,b) )_{\lambda\mu} :=
    v(r^{-1}X;r,q,p,t,ar^2,br)_{\lambda\mu}
\end{equation}
for any $v\in V$ where $r^{-1}X=\{r^{-1}x_1, r^{-1}x_2,\ldots
\}$. Then $\Omega^n$ is the algebra generated by the composite operators
  \begin{equation}
    \oW{x_i}:=\mathfrak{X_i}\mathfrak{M}
  \end{equation}
for all $i\in [n]$ and $\mathfrak{I}$, where the multiplication operation
is composition.
\end{definition}

Consider the algebra representation $\pi: \Omega^n
  \mapsto \mathrm{End} V$ of $\Omega^n$ on the space $V$
defined~by
\begin{equation}
  \pi(f)v =
f v
\end{equation}
for $f\in\Omega^n$ and $v\in V$. Then for
every linear combination (which is a polynomial due to
the symmetry proved below) $f$ in operators
$\mathfrak{W}(x_1)$, $\mathfrak{W}(x_2)$, $\ldots$, $
\mathfrak{W}(x_n)$ in $\Omega^n$, there exists an infinite lower
triangular matrix associated to it, namely $\pi(f)I = f I$. This leads
to our definition of the recursion formula for $\omega_\lambda$.
\begin{definition}
\label{definition:frakW}
With the notation as above and $z=(x_1, \ldots, x_n)$, we set
\begin{equation}
\label{frakW}
  \W{z}:= \mathfrak{W}(z^t) \cdot I
  = \W{r^{1-n} x_1} \cdots \W{r^{-1}x_{n-1}} \W{x_n}
\end{equation}
where $z^t=(x_n,x_{n-1}, \ldots, x_1)$ and $\mathfrak{W}(z^t)
:=\prod_{i=1}^n \mathfrak{W}(x_{n-i+1})$.
\end{definition}

Entrywise writing of~(\ref{frakW}) gives the recursion formula for
the \Jc\
\begin{equation}
\label{recurrence22}
  \omega_{\lambda/\tau}(y,z; r; a,b) := \sum_\mu
  \omega_{\lambda/\mu}(r^{-k}y; r; ar^{2k},
  b r^k ) \, \omega_{\mu/\tau}(z; r; a, b)
\end{equation}
where $y=(x_{1},\ldots, x_{n-k})\in\mathbb{C}^{n-k}$ and
$z=(x_{n-k+1},\ldots, x_n)\in\mathbb{C}^k$.

Using the recurrence relation~(\ref{recurrence22}) we can extend the
definition of $\omega_{\lambda/\mu}(x;a,b)$
from the single variable $x\in\mathbb{C}$ case to the multivariable
$\omega_{\lambda/\mu}(z; r; a,b)$ case with arbitrary number of variables
$z = (x_1,\ldots, x_n)\in\mathbb{C}^n$. The fact that multivariable
$\omega_{\lambda/\mu}(z; r; a,b)$ is symmetric will be proved below.

We set $\omega_{\lambda}(z; r; a,b) :=
\omega_{\lambda/0}(z; r; a,b)$. In the case of one variable
$x\in\mathbb{C}$ we get
\begin{equation}
   \omega_\lambda(x; r; a,b):=\omega_{\lambda/0}(x; r, q,p,t ;a,b) =
     \dfrac{(x^{-1})_\lambda (ax)_\lambda} {(qbx)_\lambda (q
     b /ax )_\lambda}
\end{equation}
which is independent of $r$.

The definition of the elliptic \Jc\ directly implies the following
properties.

\begin{lemma}
\label{vanishingOmega}
For $n$--part partitions $\lambda$ and $\mu$, if
$\mu\nsubseteq\lambda$, then
\begin{equation}
\omega_{\lambda/\mu}(x; r; a,b)=0.
\end{equation}
\end{lemma}

\begin{proof}
The proof follows from the vanishing
property~(\ref{FundamentalVanishing}) of $W_\lambda$ function which
says that for an $n$--part partitions $\lambda$,
$W_{\mu}(q^{\lambda}t^{\delta(n)};q, p, t, a, b) =0$ if
$\mu\nsubseteq\lambda$.
\end{proof}

We next establish the fact that the definition of
$\omega_{\lambda/\mu}$ is invariant under any
representation of $\lambda$. More precisely, we have
\begin{lemma}
\label{invariantRepLambda}
Let $\lambda= (\lambda_1,\lambda_2,\ldots, \lambda_n)$ be an $n$--part
partition with $\lambda_{n-m+1}\neq 0$ and $\lambda_{n-m+1} = \ldots =
\lambda_n = 0$ where $1\leq m\leq n$. Then we have
\begin{equation}
  \omega_{\lambda/\mu}(x; r; a,b) = \omega_{\hat{\lambda}/\hat{\mu}}(x; r;
  a,b).
\end{equation}
where $\hat{\lambda}$ and $\hat{\mu}$ denotes the $(n-m)$--part partition obtained by 
dropping the last $m$ zero parts from $\lambda$ and $\mu$ respectively.  That is,
$\hat{\lambda}=(\lambda_1,\ldots, \lambda_{n-m})$ and
$\hat{\mu}=(\mu_1,\ldots, \mu_{n-m})$. 
\end{lemma}

\begin{proof}
We first recall~(\ref{eq:ShiftedStable}) that for an $n$--part
partition $\nu= (\nu_1,\nu_2,\ldots, \nu_n)$ with  $\nu_{n-m+1} =
\ldots = \nu_n = 0$ where $1\leq m\leq n$, we have
\begin{equation}
W_{\nu} (z; q, p, t, a, b) = W_{\hat\nu} (yt^{-m};q, p, t, at^{2m}, bt^m)
\end{equation}
where $\nu'=(\nu_1,\ldots,\nu_{n-m})$ and $z =
(y,t^{\delta(m)})\in\mathbb{C}^n$ with $y\in\mathbb{C}^{n-m}$.

Assume that $\mu\subseteq\lambda$, for
    otherwise $\omega_{\lambda/\mu}$ will vanish by
    Lemma~(\ref{vanishingOmega}). Then if we set $z=q^{\lambda}
    t^{\delta(n)}$ and $y=(q^{\lambda_1}t^{n-1},\ldots,
    q^{\lambda_{n-m}}t^{m} )$ and use the property discussed above, we
    get
\begin{multline}
  W_{\mu}(q^{\lambda} t^{\delta(n)};q, p, t, bt^{2-2n}, br^{-1}t^{1-n})
  \\=  W_{\mu} ((y,t^{\delta(m)}) ;q, p, t, bt^{2-2n}, br^{-1}t^{1-n}) \\
%= W_{\mu'} (yt^{-m} ;q, p, t, bt^{2-2n+2m}, br^{-1}t^{1-n+m})
= W_{\hat\mu} (q^{\hat\lambda} t^{\delta(n-m)} ;q, p, t, bt^{2-2(n-m)},
br^{-1}t^{1-(n-m)})
\end{multline}
It remains to show that the factor
\begin{equation}
   \dfrac{(br^{-1}t^{1-n})_{\mu}} { (qt^{n-1})_{\mu}}
  \prod_{1\leq i< j \leq n} \left\{ \dfrac{ (qt^{j-i})_{\mu_i - \mu_j}
    } { (qt^{j-i-1})_{\mu_i - \mu_j} } \dfrac{
    (br^{-1}t^{3-i-j})_{\mu_i + \mu_j} } { (br^{-1}t^{2-i-j})_{\mu_i +
    \mu_j} } \right\}
\end{equation}
is also invariant under any representation of $\lambda$. It is
clear that for $k<n$ we have
\begin{equation}
\label{intoThree}
  \prod_{1\leq i< j \leq n} c_{ij} = \prod_{1\leq i< j \leq k} c_{ij}
  \cdot \prod_{1\leq i\leq k<j \leq n} c_{ij} \cdot \prod_{k< i< j
  \leq n} c_{ij}
\end{equation}
for an array of objects $c_{ij}$. If we set $k=n-m$ and
\begin{equation}
  c_{ij} := \left\{ \dfrac{
    (qt^{j-i})_{\mu_i - \mu_j} } { (qt^{j-i-1})_{\mu_i - \mu_j} } \dfrac{
    (br^{-1}t^{3-i-j})_{\mu_i + \mu_j} } { (br^{-1}t^{2-i-j})_{\mu_i +
    \mu_j} } \right\}
\end{equation}
we see that the last factor on the \rhs\ of the
product~(\ref{intoThree}) drops, because
$\mu_i - \mu_j =\mu_i + \mu_j =0$ for all $k< i< j\leq n$. For the
second factor in~(\ref{intoThree}) we have $\mu_i - \mu_j =\mu_i $,
since $\mu_j=0$ when $k<j\leq n$. Therefore, this middle factor
in the \rhs\ of~(\ref{intoThree}) becomes
\begin{equation}
   \prod_{1\leq i\leq k <j \leq n} \left\{ \dfrac{
    (qt^{j-i})_{\mu_i } } { (qt^{j-i-1})_{\mu_i } } \dfrac{
    (br^{-1}t^{3-i-j})_{\mu_i } } { (br^{-1}t^{2-i-j})_{\mu_i } }
    \right\} = \dfrac{(br^{-1}t^{1-k})_\mu
    (qt^{n-1})_\mu}{(qt^{k-1})_\mu (br^{-1}t^{1-n})_\mu }
\end{equation}
so that we have
\begin{equation}
  \prod_{1\leq i< j \leq n} c_{ij} = \prod_{1\leq i< j \leq k} c_{ij}
  \cdot \dfrac{(br^{-1}t^{1-k})_\mu
    (qt^{n-1})_\mu}{(qt^{k-1})_\mu (br^{-1}t^{1-n})_\mu }
\end{equation}
Therefore the factor
\begin{equation}
  \dfrac{(br^{-1}t^{1-n})_{\mu}} { (qt^{n-1})_{\mu}}
  \prod_{1\leq i< j \leq n} c_{ij} = \dfrac{(br^{-1}t^{1-k})_\mu }{
  (qt^{k-1})_\mu }  \prod_{1\leq i< j \leq k} c_{ij}
\end{equation}
is also invariant under any representation of $\lambda$. Other
factors of the form $(u)_\lambda$ or $(u)_\mu$ and the
very well--poised factor depend only on the length $\ell(\lambda)$
and $\ell(\mu)$ of $\lambda$ and $\mu$ respectively, and not on $n$.
\end{proof}

This shows that the apparent $n$ dependency of the
definition~(\ref{eq:omega{lambda,mu}}) is not essential in the sense
that $n$ could be any integer such that $n\geq \ell(\lambda)$. Other
properties of $\omega_{\lambda/\mu}$ will be proved after the proof of
the main theorem in the next section.

\section{$W$--\Js}
\label{section3}
We will prove in this section a number of the important properties of
the elliptic combinatorial $W$ functions.
We begin with a result that gives a special evaluation of the $W$
function, corresponding to a special evaluation of Okounkov's
$P^\star_\lambda$ functions~\cite{Okounkov1}.
\begin{prop}
\label{specialevalW}  
For any n-part partition $\lambda$
  and $a,b,q, p, t\in \mathbb{C}$ we have
\begin{equation}
\label{eq:lambda_eval}
W_\lambda(q^{\lambda}t^{\delta(n)};q,p,t,a,b) = N(\lambda,n;q,p,t,a,b)
\end{equation}
\end{prop}

\begin{proof}
Equation~(\ref{eq:lambda_eval}) is a direct consequence of repeated
applications of identities (\ref{reduction}) and
(\ref{eq:ShiftedStable}).
\end{proof}

In the proof of the main theorem in this section, we'll use a
generalization of a series transformation identity known
as \BT\ as given in~\cite{Coskun}. It 
can be stated as follows: Let $\alpha, \delta$, and $m\in V$ where $V$
is defined as in 
Definition~(\ref{Valgebra}). That is, $\alpha$,
$\delta$ and $m$ are infinite lower triangular matrices (doubly
indexed array of objects) whose rows and columns are indexed by
partitions \wrt\ the partial inclusion
ordering~(\ref{partialordering}). Denote by $\beta$ and
$\gamma$ the product matrices  $\beta = m \, \alpha$ and $\gamma =
\delta \, m$, where the matrix multiplication in $V$ is defined
by~(\ref{multiplication}). Then it is obvious that $\gamma \, \alpha =
\delta \, \beta$. Since this
identity involves three matrices on both sides, we have two ways of
summing the series due to the associativity of matrix
multiplication. This observation is stated as a change of summation
order in
\begin{lemma}[\BT]
Suppose we are given three infinite lower triangular matrices
$\alpha$, $\delta$ and $m$ in $V$.
If we define $\beta$ and $\gamma$ to be matrices with entries
\begin{equation}
\label{definitionForBT}
\beta_{\lambda\tau} = \sum_{\substack{\mu \\ \tau \subseteq
    \mu\subseteq \lambda}} m_{\lambda\mu} \, \alpha_{\mu\tau}  ,
\quad \mathrm{and} \quad
\gamma_{\nu\lambda} = \sum_{\substack{\mu \\ \lambda \subseteq
    \mu\subseteq \nu}}
\delta_{\nu\mu} \, m_{\mu\lambda}
\end{equation}
then we have
\begin{equation}
\label{BT}
\sum_{\substack{\lambda\\ \tau \subseteq
    \lambda \subseteq \nu }} \gamma_{\nu\lambda}
\, \alpha_{\lambda\tau}  =
\sum_{\substack{\lambda\\ \tau \subseteq
    \lambda \subseteq \nu } }
\delta_{\nu\lambda} \beta_{\lambda \tau}
\end{equation}
\end{lemma}
\begin{proof}
We have that
\begin{multline}
\sum_{\substack{\lambda\\ \tau \subseteq
    \lambda \subseteq \nu }} \gamma_{\nu\lambda}
\, \alpha_{\lambda\tau}  =\sum_{\substack{\lambda\\ \tau \subseteq
    \lambda \subseteq \nu }} \alpha_{\lambda\tau} \sum_{\substack{\mu
    \\ \lambda  \subseteq \mu\subseteq \nu}} \delta_{\nu\mu} \,
    m_{\mu\lambda} \,  \\
= \sum_{\substack{\mu \\ \tau \subseteq
    \mu \subseteq \nu }} \delta_{\nu\mu} \sum_{\substack{\lambda \\ \tau
    \subseteq \lambda\subseteq \mu}}
  m_{\mu\lambda} \, \alpha_{\lambda\tau}
= \sum_{\substack{\mu \\ \tau \subseteq
    \mu \subseteq \nu }} \delta_{\nu\mu} \beta_{\mu\tau}
\end{multline}
where the second step follows from the associativity of matrix
multiplication as explained above.
\end{proof}

In particular, if we set $\tau =0$ and send $\nu \rightarrow \infty$
in~(\ref{definitionForBT}) and~(\ref{BT}), we get, in the \od\ case, the
statement of the classical \BT. Because of
the limits involved, we require that certain convergence conditions are
satisfied so that $\gamma$ is well--defined and the change of
summation makes sense.

We now state one of the main theorems of the paper: a nine step
theorem proved by induction on the number of parts
(including zero parts at the end) of the partition $\lambda$.

\begin{theorem}
\label{thm:bigInduction}
In the following assume partitions $\lambda,\,\nu,\,\tau,\,\mu$ have
at most $n$ parts,
%(possibly including zero parts),
$x\in\mathbb{C}^n$ and
$q,p,t,r,s,a,a^{\prime}, b,c,d,e,f,g\in \mathbb{C}$.
\begin{enumerate}
\item The combinatorial $W_\lambda$ function satisfy conditions
(\ref{def:item1})-- (\ref{def:item4}).
\item A $W$ function symmetry identity:
\begin{multline}
\label{eq:duality}
W_{\lambda}\left(k^{-1}q^\nu t^\delta;q,,p,t,k^2a,kb\right)
\cdot\dfrac{(qbt^{n-1})_\lambda (qb/a)_\lambda}{(k)_\lambda
  (kat^{n-1})_\lambda} \\
\cdot \prod_{1\leq i < j\leq n} \left\{ \dfrac{(t^{j-i})_{\lambda_i
    -\lambda_j}(qa^{\prime}t^{2n-i-j-1})_{\lambda_i+\lambda_j}}
{(t^{j-i+1})_{\lambda_i
-\lambda_j}(qa^{\prime}t^{2n-i-j})_{\lambda_i+\lambda_j}} \right\}\\
=W_{\nu}\left(h^{-1}q^\lambda t^\delta;q, p, t,h^2a^{\prime},hb\right)
\cdot \dfrac{(qbt^{n-1})_\nu (qb/a^{\prime})_\nu}{(h)_\nu
  (ha^{\prime}t^{n-1})_\nu} \\
\cdot \prod_{1\leq i < j\leq n} \left\{ \dfrac{(t^{j-i})_{\nu_i
      -\nu_j}(qat^{2n-i-j-1})_{\nu_i+\nu_j}} {(t^{j-i+1})_{\nu_i
-\nu_j}(qat^{2n-i-j})_{\nu_i+\nu_j}} \right\}
\end{multline}
where $k=a^{\prime}t^{n-1}/b$ and $h = at^{n-1}/b$.
\item A $W$ function Jackson sum:
\begin{multline}
\label{eq:W_Jackson}
W_{\lambda}(s^{-1}x;q,p,t,at^{-2n}s^2,bt^{-n}s)\\
=\dfrac{(s)_{\lambda}(ast^{-n-1})_{\lambda}}
{(qbt^{-1})_{\lambda}(qbt^n/a)_{\lambda}}
\cdot \prod_{1\leq i < j\leq n}
\left\{\dfrac{(t^{j-i+1})_{\lambda_i-\lambda_j}
(qbst^{-i-j+1})_{\lambda_i+\lambda_j}}
{(t^{j-i})_{\lambda_i-\lambda_j} (qbst^{-i-j})_{\lambda_i +
    \lambda_j}} \right\}  \\
\cdot\sum_{\mu \subseteq \lambda}
  \dfrac{(bt^{-n})_{\mu}(qbt^n/(as))_\mu}{(qt^{n-1})_\mu
(ast^{-n-1})_\mu} \cdot \prod_{i=1}^n \left\{
  \dfrac{E(bt^{1-2i}q^{2\mu_i})}{E(bt^{1-2i})}(qt^{2i-2})^{\mu_i}
\right\} \\ \cdot \prod_{1\leq i < j \leq n}
\left\{ \dfrac{(t^{j-i})_{\mu_i -\mu_j} (qt^{j-i})_{\mu_i
      -\mu_j}}{(qt^{j-i-1})_{\mu_i -\mu_j}(t^{j-i+1})_{\mu_i
      -\mu_j}} \dfrac{(bqt^{-i-j})_{\mu_i+\mu_j}
    (bt^{-i-j+2})_{\mu_i+\mu_j}} {(bt^{-i-j+1})_{\mu_i+\mu_j}
(qbt^{-i-j+1})_{\mu_i+\mu_j}} \right\} \\
\cdot W_\mu(q^\lambda t^{\delta(n)};q,p,t,bst^{1-2n},bt^{-n})\cdot
W_\mu(x;q,p,t,at^{-2n},bt^{-n})
\end{multline}
\item A \ci\ for $\omega$ functions:
\begin{multline}
\omega_{\nu/\tau}((rs)^{-1};rs; a(rs)^2,brs)\\ =\sum_{\tau\subseteq
  \lambda\subseteq \nu} \omega_{\nu/\lambda}(s^{-1};s; a(rs)^2,brs)
\cdot  \omega_{\lambda/\tau}(r^{-1};r; ar^2,br)
\end{multline}

\item A $BC_n$ generalization of Bailey's \tns\ transformation
\begin{multline}
\dfrac{(qb/g, qb/f, q\gamma, q\gamma/fg)_{\nu}} {(qb, qb/fg, q\gamma/
     f, q\gamma /g)_{\nu} } \dfrac{(\gamma c/b,\gamma e/b)_{\tau}}{(c,
     e)_{\tau}} \\
\cdot {_{10}\Phi_9} \left[b, c,d,e,f,g;\, q\gamma b, \nu \|\, b/d,
     \tau \right] \\ = {_{10}\Phi_9} \left[\gamma, \gamma
  c/b, \gamma d/b, \gamma e/b, f, g;\,
  q\gamma b/fg, \nu \|\, b/d, \tau \right]
\end{multline}
where $\gamma = qb^2/ cde$.

\item The multivariable $\omega_{\lambda/\mu}$ function is symmetric. 
Namely,
\begin{equation}
\omega_{\lambda/\mu}(x_1,\ldots, x_k; r; a, b) =
\omega_{\lambda/\mu}(x_{\sigma(1)},\ldots, x_{\sigma(k)}; r; a, b)
\end{equation}
for any~$\sigma\in~\mathcal{S}$, the symmetric group.

\item A skew $\omega$ function special evaluation identity:
\begin{multline}
\label{omega(xr^delta)}
\omega_{\lambda/\mu}(xr^{\delta(m)};r; a,b) =
    \dfrac{(x^{-1}, axr^{m-1})_{\lambda}} {(qbr^{m-1}x,
    qb/(ax))_{\lambda} }   \dfrac{(qbr^{-1}x,
    qb/(axr^m))_{\mu} } {(x^{-1}, axr^{m-1})_{\mu} } \\  \cdot
    \dfrac{(r^m, qbr^{m-2})_{\lambda} }
    {(qbr^{-1}, r)_{\lambda} }  \dfrac{(r,
    br^{-1}t^{1-n})_{\mu}}{(qbr^{-2}, qt^{n-1})_{\mu}}
\cdot \prod_{i=1}^{n}\left\{ \dfrac{E(br^{-1}t^{2-2i} q^{2\mu_i})}
    {E(br^{-1}t^{2-2i})}  \left(qt^{2i-2} \right)^{\mu_i} \right\}\\
  \cdot\prod_{1\leq i< j\leq n}\left\{\dfrac{(qt^{j-i})_{\mu_i -\mu_j}
    } { (qt^{j-i-1})_{\mu_i - \mu_j} } \dfrac{
    (br^{-1}t^{3-i-j})_{\mu_i + \mu_j} } { (br^{-1}t^{2-i-j})_{\mu_i+
    \mu_j} } \right\} \\
\cdot W_{\mu} (q^{\lambda}t^{\delta(n)};q,p, t, br^{m-1}t^{2-2n},
    br^{-1}t^{1-n})
\end{multline}
where $x\in \mathbb{C}$ and $xr^{\delta(m)} = (xr^{m-1},\ldots, xr,
x)$ for a positive integer $m$.

\item The skew $W$ function $W_{\lambda/\mu}(x,y;q,p, t,a,b)$ is
  symmetric in the variables $x,y\in \mathbb{C}$.
\item An ``extended''
  $W$ Jackson sum: i.e. the identity (\ref{eq:W_Jackson}) for n+1 part
  partitions
$\lambda = (\lambda_1,\ldots,\lambda_{n+1})$ with the last part
$\lambda_{n+1} = 0$.
\end{enumerate}
\end{theorem}
\begin{proof}
When $\lambda = (0)$, the identities in Theorem
\ref{thm:bigInduction} are trivially true.  We now do induction on
$n$, the number of parts of $\lambda$ assuming that at least one part
is nonzero.

Observe that in step one above, identities
(\ref{def:item2})--(\ref{def:item4}) have already been proven.
Identity (\ref{def:item2}) is simply Theorem
\ref{thm:vanishingW}. Identity (\ref{def:item3}) is a consequence of
the $W$ reversal identity (\ref{eq:reversal}) and Theorem
\ref{thm:vanishingW}. Finally identity (\ref{def:item4}) is proved in
Proposition \ref{specialevalW}.  This leaves identity
(\ref{def:item1}) to be proved.  The combinatorial definition of the
$W$ function given in (\ref{definitionHfactor}), (\ref{definitionSkewW}), 
and
(\ref{eqWrecurrence}) shows that the function
$W_{\lambda}(x_1,\ldots,x_n;q,p,t,a,b)\in V_n[\lambda_1]$. To further
show that $W_{\lambda}\in V_n^{\mathcal{S}}[\lambda_1]$, we need to
demonstrate that the $W_\lambda$ function is invariant in the
variables $\{x_1,\ldots, x_n\}$ under the action of the group
$\mathcal{S}$. The difficult condition here is to show that the
$W_\lambda$ function is invariant under the transpositions
$$x_i \leftrightarrow x_j, \mbox{ for } 1\leq i< j\leq n.$$ Using
induction on n and the $y\in\mathbb{C}^2$ case of the generalized $W$
recurrence identity
\begin{multline}\label{eq:genWrecurrence}
W_{\lambda/\mu}(y,z_1,z_2,\ldots,z_\ell;q, p, t, a, b) \\ =
\sum_{\nu\subseteq \lambda} W_{\lambda/\nu}(yt^{-\ell};q, p, t,
at^{2\ell}, bt^\ell) \, W_{\nu/\mu}(z_1,\ldots, z_\ell;q, p, t, a, b),
\end{multline}
it will suffice to show that $W_{\lambda}(x_1,\ldots,x_n;q,p,t,a,b)$ is
invariant under the transposition $x_1 \leftrightarrow x_2$. When
$n=1$ there is nothing to prove, so we assume $n\geq 2$.  We can also
assume that $\lambda_n =0$, because if $\lambda_n >0$ then by the
$W_\lambda$ reduction identity (\ref{reduction}) we can reduce to the
case $\lambda_n =0$.  Using the recurrence identity
(\ref{eq:genWrecurrence}), the symmetry of the $W_\lambda$ function in
$x_1$ and $x_2$ follows from a proof of the symmetry of the function
$W_{\lambda/\nu}(x_1t^{-(n-2)}, x_2t^{-(n-2)};q, p, t, at^{2(n-1)},
bt^{n-1})$ in the variables $x_1$ and $x_2$,where  assume that
$\lambda_n =0$ and hence $\nu_n = 0$.   If $\hat{\lambda}$ and
$\hat{\nu}$ are the n-1 part partitions obtained by dropping the nth
part of $\lambda$ and $\nu$, we observe that
\begin{multline}\label{eq:SkewLastPart0}
W_{\lambda/\nu}(x_1t^{-(n-2)}, x_2t^{-(n-2)};q, p, t, at^{2(n-2)},
bt^{n-2})\\ = W_{\hat{\lambda}/\hat{\nu}}(x_1t^{-(n-2)},
x_2t^{-(n-2)};q, p, t, at^{2(n-2)}, bt^{n-2}).
\end{multline}
Applying the step seven (n-1 part partitions) inductive assumption to
the right-hand-side of (\ref{eq:SkewLastPart0}), we finish the proof
that $W_{\lambda}(x_1,\ldots,x_n;q,.pt,a,b)\in
V_n^{\mathcal{S}}[\lambda_1]$ and step one of the n part partitions
case of Theorem \ref{thm:bigInduction}.

To prove step two, identity (\ref{eq:duality}), we will do a double
induction, first on the last part $\lambda_n$ of $\lambda$, and
finally on the last part $\nu_n$ of $\nu$. By the general inductive
assumption for Theorem \ref{thm:bigInduction}, we also assume that
identity (\ref{eq:duality}) is true when $\lambda$ and $\nu$ have only
n-1 parts. The general inductive assumption also implies that the
``extended'' $W$ Jackson sum (\ref{eq:W_Jackson}) is true for n-part
partitions $\lambda$ with last part zero, $\lambda_n=0$.  Substituting
$at^{2n}$ in place of $a$ and $bt^n$ in place of $b$, and setting
$x=q^{\nu}t^{\delta(n)}$, $s=k$, $k=a^{\prime}t^{n-1}/b$ and $h =
at^{n-1}/b$ in this ``extended'' $W$ Jackson sum identity, the
right-hand-side of the $W$ Jackson sum (\ref{eq:W_Jackson}) exhibits a
symmetry between $\lambda$, $a$, $k$ and $\nu$, $a^{\prime}$,
$h$. This proves the n-part partition identity (\ref{eq:duality}) when
$\lambda_n=\nu_n=0$.  If $\lambda_n > 0$, then using identity
(\ref{reduction}) the left-hand side of (\ref{eq:duality}) becomes
\begin{multline}\label{lsym}
\prod_{i=1}^n\frac{E(kq^{-\nu_i}t^{i-n})
E(kaq^{\nu_i}t^{n-i})E(qbt^{n-i})E(qba^{-1}t^{1-i})}
{E(kt^{1-i})E(kat^{n-i})E(qbq^{\nu_i}t^{n-i})
E(qba^{-1}q^{-\nu_i}t^{i-n})}\\ \cdot\prod_{1\leq i<j\leq n}
\frac{(t^{j-i})_{\lambda_i -\lambda_j}
(q^3a^{\prime}t^{2n-i-j-1})_{\lambda_i+\lambda_j-2}}
{(t^{j-i+1})_{\lambda_i-\lambda_j}(q^3a^{\prime}
t^{2n-i-j})_{\lambda_i+\lambda_j-2}}\\ \cdot
\frac{(q^2bt^{n-1})_{\lambda-1^n} (q^2b/a)_{\lambda-1^n}}
{(qk)_{\lambda-1^n} (qkat^{n-1})_{\lambda-1^n}}
W_{\lambda-1^n}\left(k^{-1}q^\nu t^\delta
q^{-1};q, p, t,k^2aq^2,kbq^2\right).
\end{multline}
where $\lambda-1^n$ is the partition
$(\lambda_1-1,\ldots,\lambda_n-1)$ and we have simplified by using the
identity
\begin{multline}
\frac{(qkbt^{j-2i})_{2}(qa^{\prime}t^{2n-i-j-1})_2}
{(qkbt^{j-1-2i})_{2}(qa^{\prime}t^{2n-i-j})_2}=
\frac{(qa^{\prime}t^{n+j-2i-1})_{2}(qa^{\prime}t^{2n-i-j-1})_2}
{(qa^{\prime}t^{n+j-2-2i})_{2}(qa^{\prime}t^{2n-i-j})_2}\\
=\frac{(qa^{\prime}t^{n-i-1};t)_{n-i}(q^2a^{\prime}t^{n-i-1};t)_{n-i}}
{(qa^{\prime}t^{n-i-1};t)_{n-i}(q^2a^{\prime}t^{n-i-1};t)_{n-i}}=1.
\end{multline}
By induction on $\lambda_n$, we can evaluate the $W_{\lambda-1^n}$
function in (\ref{lsym}) from the ${\lambda-1^n}$ case of
(\ref{eq:duality}) with $b,k,h,a^{\prime}$ replaced by $qb,qk,
q^{-1}h, q^2a^{\prime}$ respectively and $\nu$ and $a$ unchanged. We
obtain
\begin{multline}\label{lsym2}
W_{k^{-1}q^\nu t^\delta
q^{-1};q,p, \lambda-1^n}\left(t,k^2aq^2,kbq^2\right)\\ = 
\frac{(qk)_{\lambda-1^n}
(qkat^{n-1})_{\lambda-1^n}}{(q^2bt^{n-1})_{\lambda-1^n}
(q^2b/a)_{\lambda-1^n}} \prod_{1\leq i < j\leq
n}\frac{(t^{j-i+1})_{\lambda_i
-\lambda_j}(q^3a^{\prime}t^{2n-i-j})_{\lambda_i+\lambda_j-2}}
{(t^{j-i})_{\lambda_i
-\lambda_j}(q^3a^{\prime}t^{2n-i-j-1})_{\lambda_i+\lambda_j-2}}\\
\cdot\frac{(q^2bt^{n-1})_\nu (b/a^{\prime})_\nu} {(q^{-1}h)_\nu
(qha^{\prime}t^{n-1})_\nu} \prod_{1\leq i < j\leq
n}\frac{(t^{j-i})_{\nu_i -\nu_j}(qat^{2n-i-j-1})_{\nu_i+\nu_j}}
{(t^{j-i+1})_{\nu_i -\nu_j}(qat^{2n-i-j})_{\nu_i+\nu_j}}\\ \cdot
W_{\nu}\left(h^{-1}q^\lambda t^\delta;q,p, t,h^2a^{\prime},hb;\right).
\end{multline}
Substituting identity (\ref{lsym2}) into the expression (\ref{lsym})
and simplifying further, we prove identity (\ref{eq:duality}) in the
$\lambda_n >0$ case.

The case with $\nu_n > 0$ can be proven in exactly the same way.
Therefore we reduce to the previously proven case $\lambda_n = \nu_n=
0$.  This completes the proof of step two of Theorem
\ref{thm:bigInduction}.

We now prove induction step three of Theorem \ref{thm:bigInduction},
the $W$ function Jackson sum (\ref{eq:W_Jackson}) in the case when the
variable $x\in \mathbb{C}^n$, treating the partition $\lambda$ as an
n-tuple.  We will do a further induction on $\lambda$ using reverse
lexicographical order. We will use the notation $\mu<\lambda$ to
denote $\mu$ is less than $\lambda$ in reverse lexicographical order
(and similarly $\mu\leq \lambda$ also includes the case
$\mu=\lambda$).  The inductive assumption is that identity
(\ref{eq:W_Jackson}) is true whenever $\lambda<\gamma$ for some
partition $\gamma$ with $l(\gamma)\leqq n$. We now consider the case
$\lambda =\gamma$ and proceed to show that identity
(\ref{eq:W_Jackson}) is true in this case as well. Note that we can
begin our induction in the trivial case of the partition (0).  We
begin with the following
\begin{lemma}\label{lem:Wexpans}
  \begin{equation}
    \label{eq:Wexpans}
    W_\lambda (s^{-1}x; q,,p,t,at^{-2n}s^2,bs) =
\sum_{\mu\leq\lambda}c_\mu W_\mu(x;q,p,t,at^{-2n},b),
  \end{equation}
where $x\in \mathbb{C}^n$, $s\in \mathbb{C}$, $c_\mu \in \mathbb{F}$,
and the sum is over partitions $\mu$ with $l(\mu)\leq n$ and $
\mu\leq\lambda$.
\end{lemma}
\begin{proof}
  Since $W_\lambda (s^{-1}x; q,p,t,at^{-2n}s^2,bt^{-n}s) \in
V_n^{\mathcal{S}}[\lambda_1]$, it follows that
\begin{equation}
  \label{eq:first_part_expand}
  W_\lambda (s^{-1}x; q,p,t,at^{-2n}s^2,bt^{-n}s) = \sum_\mu c_\mu
W_\mu(x;q,p,t,at^{-2n},bt^{-n}),
\end{equation}
where the sum is over partitions $\mu$ with $l(\mu)\leq n$ and
$\mu_1\leq \lambda_1$.  We will denote the residue of any function
$f\in V_n^{\mathcal{S}}[\lambda_1]$ at the pole $x_n =
(bt^{-n}q^{\lambda_1})^{-1}$ in the variable $x_n$ by $
R_{n,\lambda_1}(f)$ and consider the residues of both sides of
equation (\ref{eq:first_part_expand}).  First observe that
$R_{n,\lambda_1}$ induces an injection
\begin{equation}
  \label{eq:res_inj}
  R_{n,\lambda_1}:
(V_n^{\mathcal{S}}[\lambda_1]/V_n^{\mathcal{S}}[\lambda_1-1])
\longrightarrow V_{n-1}^{\mathcal{S}}[\lambda_1].
\end{equation}
Next observe that
\begin{multline}
  \label{def:W_tilde}
\tilde{W}_\lambda ^{(n-1)} = I_{n,\lambda_1} (W_\lambda (s^{-1}x;
q,p,t,at^{-2n}s^2,bt^{-n}s))\\
=\prod_{i=1}^{n-1}\dfrac{(qbt^{-n}x_i)_{\lambda_1}
(qbt^{-n}/(at^{-2n}x_i)_{\lambda_1}}
{(qbt^{-n}x_i/t)_{\lambda_1}(qbt^{-n}/(at^{-2n}x_it)_{\lambda_1}}\\
\cdot R_{n,\lambda_1}(W_\lambda (s^{-1}x; q, p, t,at^{-2n}s^2,bt^{-n}s))
\in V_{n-1}^{\mathcal{S}}[\lambda_2].
\end{multline}
This implies that on the right-hand side of the expansion
(\ref{eq:first_part_expand}), the coefficients $c_\mu = 0$ unless
$(\mu_1<\lambda_1)$ or $(\mu_1 = \lambda_1\mbox{ and
}\mu_2\leq\lambda_2)$. Observe similarly that
\begin{multline}\label{def:Imap2}
\tilde{W}_\lambda ^{(n-2)} =
I_{n-1,\lambda_2}(\tilde{W}_\lambda^{(n-1)})\\
=\prod_{i=1}^{n-2}\dfrac{(qbt^{-n}x_i)_{\lambda_2}
(qbt^{-n}/(at^{-2n}x_i)_{\lambda_2}}
{(qbt^{-n}x_i/t)_{\lambda_2}(qbt^{-n}/(at^{-2n}x_it)_{\lambda_2}}
\cdot R_{n-1,\lambda_2}(\tilde{W}_\lambda^{(n-1)}) \in
V_{n-2}^{\mathcal{S}}[\lambda_3],
\end{multline}
which implies that, in expansion (\ref{eq:first_part_expand}), the
coefficients $c_\mu = 0$ unless $(\mu_1<\lambda_1)$ or $(\mu_1 =
\lambda_1\mbox{ and }\mu_2<\lambda_2)$ or $(\mu_1 = \lambda_1\mbox{,
}\mu_2=\lambda_2\mbox{ and } \mu_3\leq\lambda_3).$ Iterating this
procedure then proves expansion (\ref{eq:Wexpans}) and completes the
proof of Lemma \ref{lem:Wexpans}.
\end{proof}
Using Lemma \ref{lem:Wexpans}, the vanishing theorem (Theorem
\ref{thm:vanishingW}), and identity (\ref{eq:duality}), we can
inductively compute the coefficients $c_\mu$ in expansion the
(\ref{eq:Wexpans}) for
$W_{\lambda}(s^{-1}x;q,p,t,at^{-2n}s^2,bt^{-n}s)$.  Let $\nu$ be a
partition with $l(\nu)\leq n$ and $\nu<\lambda$.  Suppose we have for
$\mu<\nu$, that all the coefficients $c_\mu$ are equal to the
corresponding coefficients of $W_\mu(x;q,p,t,at^{-2n},bt^{-n})$ on the
right-hand side of expansion (\ref{eq:W_Jackson}).  Setting $x=q^\nu
t^{\delta(n)}$ in identity (\ref{eq:W_Jackson}) will allow us to
compute $c_\nu$. By Theorem \ref{thm:vanishingW} all the $W_\mu(x)$
terms on the right-hand side of (\ref{eq:W_Jackson}) will vanish when
$\nu < \mu$.  Then the coefficient $c_\nu$ in the expansion for
$W_{\lambda}(s^{-1}x;q,p,t,at^{-2n}s^2,bt^{-n}s)$ can be computed by the
identity (\ref{eq:duality}) and the induction hypothesis for
(\ref{eq:W_Jackson}).  Identity (\ref{eq:duality}) allows us to
compute the expansion of
$W_{\lambda}(s^{-1}x;q,p,t,at^{-2n}s^2,bt^{-n}s)$ at $x=q^\nu
t^{\delta(n)}$ by means of the expansion of
$W_{\nu}\left(h^{-1}x;q,p,t,h^2asbt^{1-n},hb\right)$ at $x=q^\lambda
t^{\delta(n)}$ where $h=at^{n-1}/b$.  Since all the non-vanishing
terms in $W_\lambda$ at $x=q^\nu t^{\delta(n)}$ are known except for
the $c_\nu$ term and all the terms of the expansion
$W_{\nu}\left(h^{-1}x;q,p,t,h^2asbt^{1-n},hb\right)$ are known by the
induction hypothesis for (\ref{eq:W_Jackson}), then we can solve for
$c_\nu$.  Given the symmetry of the expansion of
$W_{\lambda}(s^{-1}x;q,p,t,at^{-2n}s^2,bt^{-n}s)$ at $x=q^\nu
t^{\delta(n)}$, it is easily seen that the coefficient $c_\nu$ agrees
with the corresponding coefficient in the $W$-Jackson sum
(\ref{eq:W_Jackson}).  By induction on the partitions $\nu$ with
$l(\nu)\leq n$ and $\nu<\lambda$, it follows that all the coefficients
$c_\nu$, $\nu<\lambda$, in the expansion of
$W_{\lambda}(s^{-1}x;q,p,t,at^{-2n}s^2,bt^{-n}s)$ agree with the
corresponding coefficients on the right-hand side of
(\ref{eq:W_Jackson}).  We are left to show that the leading
coefficient $c_\lambda$ equals the coefficient of
$W_\lambda(x;q,p,t,at^{-2n},bt^{-n})$ in (\ref{eq:W_Jackson}).

The last part of the proof of step three of Theorem
\ref{thm:bigInduction} is a direct computation of the coefficient
$c_\lambda$. This can be done by computing the result
$\tilde{W}_\lambda^(0)$ of composing the n $I$ maps defined in
(\ref{def:W_tilde}), (\ref{def:Imap2}), etc., acting on
$W_{\lambda}(s^{-1}x;q,p,t,at^{-2n}s^2,bt^{-n}s)$.  In order to
calculate the residues $R_{n-k,\lambda_{k+1}}$ for $k=0,\ldots, n-1$
we apply the combinatorial definition of $W_{\lambda}(s^{-1}x)$ given
in definitions (\ref{definitionHfactor}) --(\ref{eqWrecurrence}) and
repeatedly apply the one-dimensional Jackson sum:
\begin{multline}
  \label{eq:1dJackson}
\dfrac{(azs,z^{-1}s,qb,qb/a)_m}{(qb/(az),qbz,as,s)_m} =\sum_{k=0}^m
\dfrac{E(q^{2k}b)(az,z^{-1}b,
qb/(as),bsq^m,q^{-m})_k}{E(b)(qb/(az),qbz,q,as,q^{1-m}s^{-1},
bq^{1+m})_k}q^k,
\end{multline}
for any non-negative integer m.  The result is, after some
simplification,
\begin{multline}
  \label{leadingcoeffmap}
\tilde{W}_\lambda^(0) =
s^{|\lambda|}\dfrac{(bt^{-n},qbt^n/(as))_\lambda}
{(bt^{-n}s,qbt^n/a)_\lambda}
\cdot\prod_{i=1}^n\dfrac{(bst^{1-2i})_{2\lambda_i}}
{(bt^{1-2i})_{2\lambda_i}}\\
\cdot\prod_{1\leq i< j\leq n}\dfrac{(qbst^{1-i-j},qbt^{-i-j},
bt^{2-i-j},bst^{1-i-j})_{(\lambda_i +\lambda_j)}}
{(qbt^{1-i-j},qbst^{-i-j},bst^{2-i-j},
bt^{1-i-j})_{(\lambda_i +\lambda_j)}}.
\end{multline}
We obtain the identical result from
composing the n $I$ maps acting on the
right-hand side of (\ref{eq:W_Jackson}),
using the identity (\ref{eq:lambda_eval}) applied to the  factor
$W_\lambda(q^\lambda t^{\delta})$ as well as the
definitions (\ref{definitionHfactor})--(\ref{eqWrecurrence}).
This proves that $c_\lambda$ equals
the coefficient of $W_\lambda(x)$ on the
right-hand side of (\ref{eq:W_Jackson})
and completes the proof of step three of
Theorem \ref{thm:bigInduction}.

We are now ready to prove the induction step four
, the \ci~(\ref{CocycleIdentity}), mentioned
above. We first note that the \Js~(\ref{eq:W_Jackson}) can be
written as
\begin{equation}
\label{stillJackson}
  W_\lambda(r^{-1}z; q,p, t, a r^2, b r)
= \sum_{\mu\subseteq\lambda} c_{\lambda/\mu}(r, a ,b)
\, W_\mu(z; q, p, t, a, b)
\end{equation}
where
\begin{multline}
c_{\lambda/\mu}(r, a ,b)
:= \dfrac{[r]_\lambda [a rt^{n-1}]_\lambda}
{[qb t^{n-1}]_\lambda [qb /a]_\lambda} \frac{[b]_\mu [qb/a r]_\mu} 
{[qt^{n-1}]_\mu
  [art^{n-1}]_\mu}\, W_\mu(q^\lambda t^{\delta(n)}; q, p,t, b rt^{1-n},
b) \\ \cdot \prod_{1\leq i<
  j\leq n} \frac{(t^{j-i+1})_{\lambda_i-\lambda_j} (qb
  rt^{1+n-i-j})_{\lambda_i+\lambda_j}} {(t^{j-i})_{\lambda_i-\lambda_j}
  (qb r t^{n-i-j})_{\lambda_i+\lambda_j}} \cdot\prod_{i=1}^n
\left\{\frac{E(b t^{1+n-2i}q^{2\mu_i})}
  {E(b t^{1+n-2i})} (qt^{2i-2})^{\mu_i}\right\}  \\ \cdot \prod_{1\leq i<j 
\leq
n} \left\{\frac{(t^{j-i})_{\mu_i-\mu_j} (qt^{j-i})_{\mu_i-\mu_j}
(b qt^{n-i-j})_{\mu_i+\mu_j} (bt^{2+n-i-j})_{\mu_i+\mu_j}}
{(qt^{j-i-1})_{\mu_i-\mu_j} (t^{j-i+1})_{\mu_i-\mu_j}
(b t^{1+n-i-j})_{\mu_i+\mu_j} (b qt^{1+n-i-j})_{\mu_i+\mu_j}}
\right\}
\end{multline}
for any $r\in\mathbb{C}$. We note that
the identity~(\ref{stillJackson}) could be iterated repeatedly.
A double iteration gives
\begin{lemma}
\label{lem:extendedJsMatrixVersion1}
With the notation as above and $u,v\in\mathbb{C}$, we have
\begin{multline}
\label{extendedJsMatrixVersion1}
  c_{\nu/\mu}((uv)^{-1}, a(uv)^2, buv)  =\sum_{\mu\subseteq
  \lambda \subseteq \nu} c_{\nu/\lambda}(v^{-1},a(vu)^2,bvu) \,
  c_{\lambda/\mu}(u^{-1},au^2,bu)
\end{multline}
\end{lemma}

\begin{proof}
We start with the \Js~(\ref{stillJackson}) in the form
\begin{multline}
W_\lambda((uv)^{-1}x; q,p, t, a(uv)^2, buv) \\ = \sum_\mu
  c_{\lambda/\mu}(v^{-1}, a(uv)^2, buv) \,
   W_\mu(u^{-1}x; q, p, t, au^2, bu)
\end{multline}
Expand the $W_\lambda$ functions on both sides using the
\Js~(\ref{stillJackson}) again and write
\begin{multline}
  \sum_\tau
  c_{\lambda/\tau}((uv)^{-1},a(uv)^2,buv) \,
  W_\tau(x;q,p, t, a,b) \\ = \sum_\mu
  c_{\lambda/\mu}(v^{-1}, a(uv)^2, buv) \,
  \sum_\tau c_{\mu/\tau}(u^{-1}, au^2,bu) \, W_\tau(x;q,p, t,a,b)
\end{multline}
Switch the order of summation and compare the coefficients of the
basis functions $W_\tau(x;a,b)$ to get the identity to be proved.
\end{proof}

Note that some factors in the
expansion~(\ref{extendedJsMatrixVersion1}) cancel. We now rewrite the
above theorem~(\ref{lem:extendedJsMatrixVersion1}) in terms of the
$\omega_{\lambda/\mu}$ function in a way involving the essential
factors and independent of different representations of $\nu$ in
\begin{cor}[Cocycle Identity]
\label{cor:CocycleIdentity}
Let $\lambda$ be a partition. With the notation as above and
$u,v\in\mathbb{C}$, we have
\begin{multline}
\label{CocycleIdentity}
  \omega_{\nu/\mu}((uv)^{-1};uv;a(uv)^2, buv) \\ = \sum_{\mu\subseteq
  \lambda \subseteq \nu} \omega_{\nu/\lambda}(v^{-1};v;a(vu)^2,bvu) \,
  \omega_{\lambda/\mu}(u^{-1};u;au^2,bu)
\end{multline}
\end{cor}

\begin{proof}
It suffices to note that
\begin{equation}
  \omega_{\lambda/\mu}(u^{-1};u;au^2,bu)  = P_\lambda (b) \,
  c_{\lambda/\mu}(u, au^{-2}t^{1-n}, bu^{-1}t^{1-n}) \,
  P_\mu^{-1}(bu^{-1})
\end{equation}
where
\begin{equation}
P_\lambda (b):= \prod_{1\leq i<
  j\leq n} \frac{(t^{j-i+1})_{\lambda_i-\lambda_j} (qb
  t^{2-i-j})_{\lambda_i+\lambda_j}} {(t^{j-i})_{\lambda_i-\lambda_j}
  (qbt^{1-i-j})_{\lambda_i+\lambda_j}}
\end{equation}
The diagonal factor $P_\lambda(b)$ cancels in the expansion, and we get
the identity to be proved
after a simple reparametrization.
\end{proof}

\begin{remark}
An important application of the \ci\ is the proof of a $BC_n$ \tns\
transformation which we now obtain in step six. The \tns\
transformation is an importantant application of the
\BT~(\ref{BT}). In fact, the 
\tns\ transformation formula it is referred to as the \BT\ in the
literature.
\end{remark}

The \tns\ transformation formula is proved by summing both inner series
above in the \BT~(\ref{BT}) (i.e. computing $\beta$ and
$\gamma$) by means of the \Js. We use the \ci\ for $\omega_\lambda$
to prove a more general version of the classical \tns\ transformation in
\begin{lemma}
\label{lem:BCn10phi9t}
For any partitions $\nu$ and $\tau$
and complex parameters $b,c,d,e,f$ and $g$ we have
\begin{multline}
\label{BCn10phi9t}
\dfrac{( qb/g, qb/f, q\gamma, q\gamma/fg)_{\nu} }
     {(qb, qb/fg, q\gamma/ f, q \gamma /g)_{\nu} }  \dfrac{ (\gamma
     c/b, \gamma e/b)_{\tau}} {(c, e)_{\tau}} \\
\cdot {_{10}\Phi^n_9} \left[b, c,d,e,f,g;\, q\gamma b, \nu \|\, b/d,
     \tau \right]_{q,p,t} \\ = {_{10}\Phi^n_9} \left[\gamma, \gamma
  c/b, \gamma d/b, \gamma e/b, f, g;\,
  q\gamma b/fg, \nu \|\, b/d, \tau \right]_{q,p,t}
\end{multline}
where $\gamma = qb^2/ cde$, $\tau$ and $\nu$ are partitions of length
at most $n$ such that $\tau\subseteq \nu$.
\end{lemma}
\begin{proof}
Recall that the \ci\
for $\omega_{\lambda/\mu}$ function may be written in the form
\begin{equation}
\label{version1}
  \omega_{\lambda/\tau}((sr)^{-1};sr;as^2, bs)
=\sum_\mu
  \omega_{\lambda/\mu}(s^{-1};s;as^2,bs) \,
  \omega_{\mu/\tau}(r^{-1};r;a,b)
\end{equation}
We will use this identity to compute the inner sum (say
$\beta_\lambda$) in one side of the \BT~(\ref{BT}) and the
reparametrized version
\begin{equation}
\label{version2}
  \omega_{\nu/\lambda}((us)^{-1};us;a'u^2, bsu)
=\sum_\mu
  \omega_{\nu/\mu}(u^{-1};u;a'u^2,bsu) \,
  \omega_{\mu/\lambda}(s^{-1};s;a',bs)
\end{equation}
of this identity to compute the inner sum (say $\gamma_\lambda$) on
the other side of~(\ref{BT}).

Now we set
\begin{equation}
m_{\lambda\mu} = W_\mu(q^\lambda t^{\delta(n)}; q,p,t, bst^{2-2n},
bt^{1-n})
\end{equation}
and
\begin{multline}
\alpha_{\mu\tau}= \dfrac{ (qb/as, bt^{1-n}, r, ar^{-1})_{\mu}}{
    (as, qt^{n-1}, qbr^{-1}, qbr/a)_{\mu}} \, W_{\tau}
    (q^{\mu}t^{\delta(n)}; q ,p,t, bt^{2-2n}, br^{-1}t^{1-n}) \\  \cdot
  \prod_{i=1}^{n}\left\{ \dfrac{E(bt^{2-2i} q^{2\mu_i})}
    {E(bt^{2-2i})}  \left(qt^{2i-2} \right)^{\mu_i} \right\}
  \prod_{1\leq i< j \leq n} \left\{ \dfrac{ (qt^{j-i})_{\mu_i - \mu_j}
    } { (qt^{j-i-1})_{\mu_i - \mu_j} } \dfrac{
    (bt^{3-i-j})_{\mu_i + \mu_j} } { (bt^{2-i-j})_{\mu_i +
    \mu_j} } \right\}
\end{multline}
and also set
\begin{multline}
\delta_{\nu\mu} = \dfrac{ (qbs/a'u, bst^{1-n}, s, a's^{-1})_{\mu} }{
    (a'u, qt^{n-1}, qb, qbs^2/a']_{\mu}} \, W_{\mu}
    (q^{\nu}t^{\delta(n)};q ,p, t, bsut^{2-2n}, bst^{1-n}) \\ \cdot
  \prod_{i=1}^{n}\left\{ \dfrac{E(bst^{2-2i} q^{2\mu_i})}
    {E(bst^{2-2i})}  \left(qt^{2i-2} \right)^{\mu_i} \right\}
  \prod_{1\leq i< j \leq n} \left\{ \dfrac{ (qt^{j-i})_{\mu_i - \mu_j}
    } { (qt^{j-i-1})_{\mu_i - \mu_j} } \dfrac{
    (bst^{3-i-j})_{\mu_i + \mu_j} } { (bst^{2-i-j})_{\mu_i +
    \mu_j} } \right\}
\end{multline}
Using the first version~(\ref{version1}) of the
\ci\
we compute
\begin{multline}
\beta_{\lambda \tau}
= \dfrac{(rs, asr^{-1}, qb, qb/a)_{\lambda} } {( qbr^{-1}, qbr/a, s,
  as)_{\lambda}  } \dfrac{(qb/as, ar^{-1})_{\tau}} {(asr^{-1},
  qb/a)_{\tau}} \\   \cdot
W_{\tau} (q^{\lambda}t^{\delta(n)};q, p, t, bst^{2-2n}, br^{-1}t^{1-n})
\end{multline}
and by means of the reparametrized version~(\ref{version2}) we compute
$\gamma_\lambda$, which becomes
\begin{multline}
\gamma_{\nu\lambda}
= \dfrac{(su, a'us^{-1}, qbs, qbs/a')_{\nu} } {(qb, qbs^2/a', u,
  a'u)_{\nu}} \dfrac{(qbs/a'u , a's^{-1})_{\lambda} }{ (a'us^{-1},
  qbs/a')_{\lambda} }
\\  \cdot W_{\lambda} (q^{\nu}t^{\delta(n)};q, p, t, bsut^{2-2n}, bt^{1-n})
\end{multline}
So, by using the \BT~(\ref{BT}) we get
\begin{multline}
\label{BCn10phi9tExplicit}
\dfrac{(su, a'us^{-1}, qbs, qbs/a')_{\nu} } {(qb, qbs^2/a', u,
    a'u)_{\nu}} \dfrac{ (asr^{-1}, qb/a)_{\tau} }{ (qb/as,
    ar^{-1})_{\tau} } \\ \cdot
  \sum_{\substack{\lambda \\ \mu \subseteq
    \lambda \subseteq \nu }} \dfrac{ (bt^{1-n}, qb/as, r, ar^{-1},
    qbs/a'u, a's^{-1})_{\lambda} } { (qt^{n-1}, as, qbr^{-1}, qbr/a,
    a'us^{-1}, qbs/a')_{\lambda} }  \\ \cdot
  \prod_{i=1}^{n}\left\{ \dfrac{E(bt^{2-2i} q^{2\lambda_i})}
    {E(bt^{2-2i})}  \left(qt^{2i-2} \right)^{\lambda_i} \right\}
  \prod_{1\leq i< j \leq n} \left\{ \dfrac{ (qt^{j-i})_{\lambda_i - 
\lambda_j}
    } { (qt^{j-i-1})_{\lambda_i - \lambda_j} } \dfrac{
    (bt^{3-i-j})_{\lambda_i + \lambda_j} } { (bt^{2-i-j})_{\lambda_i +
    \lambda_j} } \right\} \\ \cdot W_{\lambda} (q^{\nu}t^{\delta(n)};
    q, p, t, bsut^{2-2n}, bt^{1-n}) \, W_{\tau}
    (q^{\lambda}t^{\delta(n)};q, p, t, bt^{2-2n}, br^{-1}t^{1-n}) \\
= \sum_{\substack{\lambda\\ \mu \subseteq
    \lambda \subseteq \nu }} \dfrac{(bst^{1-n}, rs, asr^{-1}, qb/a,
    qbs/a'u, a's^{-1})_{\lambda} }  {(qt^{n-1}, qbr^{-1} , qbr/a, as, a'u,
     qbs^2/a')_{\lambda} } \\ \cdot
\prod_{i=1}^{n}\left\{ \dfrac{E(bst^{2-2i} q^{2\lambda_i})}
    {E(bst^{2-2i})}  \left(qt^{2i-2} \right)^{\lambda_i} \right\}
  \prod_{1\leq i< j \leq n} \left\{ \dfrac{ (qt^{j-i})_{\lambda_i - 
\lambda_j}
    } { (qt^{j-i-1})_{\lambda_i - \lambda_j} } \dfrac{
    (bst^{3-i-j})_{\lambda_i + \lambda_j} } { (bst^{2-i-j})_{\lambda_i +
    \lambda_j} } \right\} \\ \cdot W_{\lambda}
    (q^{\nu}t^{\delta(n)};q, p,t, bsut^{2-2n}, bst^{1-n}) \, W_{\tau} (
    q^{\lambda}t^{\delta(n)}; q, p,t, bst^{2-2n}, br^{-1}t^{1-n})
\end{multline}
A simple reparametrization $b=b$, $c=qb/as$, $d=r$, $e=ar^{-1}$,
$f=qbs/a'u$, $g=a's^{-1}$ in the above
identity~(\ref{BCn10phi9tExplicit}) completes the proof. Using the
definition~(\ref{MBHS}) we write this result in the form given above.
\end{proof}

The identity~(\ref{BCn10phi9t}) or~(\ref{BCn10phi9tExplicit})
reduces to a $BC_n$ generalization of the Bailey's classical \tns\
transformation when we set $\tau=0$.

In step seven we prove the symmetry using the
elliptic \tns\ transformation formula. Note that the
$\omega_\lambda(z; r; a,b)$, by definition, has also the built--in
symmetry $x_i\mapsto (ax_i)^{-1}$ for each $i\in
[k]$. Therefore we have
\begin{lemma}
\label{symmetry}
Let $\lambda,\mu$ be partitions and
$z=(x_1,\ldots,x_k)\in\mathbb{C}^k$. Then
\begin{equation}
\omega_{\lambda/\mu}(x_1,\ldots, x_k; r; a, b) =
\omega_{\lambda/\mu}(x_{\sigma(1)},\ldots, x_{\sigma(k)}; r; a, b)
\end{equation}
for any~$\sigma\in~\mathcal{S}$, the \hgr, where the inversions are as
defined above.
\end{lemma}
\begin{proof}
First we show that, for $x_1, x_2\in\mathbb{C}$ we have
\begin{equation}
\omega_{\lambda/\mu}(x_1,x_2; r; a,b) = \omega_{\lambda/\mu}(x_2,x_1;
r; a,b)
\end{equation}
which follows from a double application
of the elliptic  $BC_n$ generalization of Bailey's \tns\
transformation~(\ref{BCn10phi9t}). Namely, we first set $b=b,
c=x_2^{-1}, d=ax_2, e=r, f=qbx_1/r$, 
and $g=qb/arx_1$ so that $\gamma_1 = qb^2/ ar$, 
and apply the transformation~(\ref{BCn10phi9t}) above. Then apply it
again by setting 
$b=\gamma_1, c=qb/a, d=qbx_1/r, e=qb/arx_1, f=qbx_2/r, 
g=qb/arx_2 $ this time where we compute $\gamma_2=b$. This gives the
symmetry for the two variable case 
mentioned above, and implies the commutativity of the operators  
$\mathfrak{W}(x_1)$ and $\mathfrak{W}(x_2)$.

The rest of the proof follows by induction. Assuming that the symmetry
holds for $m$ variables where $m\leq k$, we show that it also holds
for $k+1$ variables. For $z=(x_1,\ldots, x_k)\in\mathbb{C}^k$ and
$x_{k+1}\in\mathbb{C}$ we have
\begin{multline}
\W{z,x_{k+1}}_{\lambda\mu} = (\mathfrak{W}(x_{k+1},z^t) \cdot
I)_{\lambda\mu} \\
= (\W{r^{-2} x_1,\ldots,r^{-2} x_{k-1}} \W{x_n,x_{k+1}})_{\lambda\mu}\\
= (\W{r^{-2} x_1,\ldots,r^{-2} x_{k-1}} \W{x_{k+1},x_n})_{\lambda\mu}\\
= (\W{r^{-1} x_1,\ldots,r^{-1} x_{k-1},r^{-1}x_{k+1}} \W{x_n})_{\lambda\mu}
\end{multline}
by the definition of the multiplication rule and the first
part of the proof. But $\W{r^{-1} x_1,\ldots,r^{-1} x_{k-1},r^{-1}x_{k+1}}$
is symmetric by hypothesis which concludes the induction
step. This is because of the fact that the
symmetric group $S_{k+1}$ in $(k+1)$ letters is generated by the
permutation $(12\ldots k)$ and the transposition $(k
(k+1))$. Therefore, $\W{z} = \W{\sigma(z)}$ for any $\sigma\in S_k$ as
required.
\end{proof}

We now employ the \ci~(\ref{CocycleIdentity}) to
compute a closed expression (\ref{omega(xr^delta)}), step seven, for 
$\omega_{\lambda/\mu}(xr^{\delta(m)};r;a,b)$.  We get
\begin{lemma}
\label{lem:omega(xr^delta)}
Let $\lambda$ and $\mu$ be partitions,
$m\in\mathbb{Z}_{>}$, and $q, p, t, a,b,r,x\in\mathbb{C}$ be complex
parameters. Then we have
\begin{multline}
\omega_{\lambda/\mu}(xr^{\delta(m)};r,a,b) =
    \dfrac{(x^{-1}, axr^{m-1})_{\lambda}} {(qbr^{m-1}x,
    qb/(ax))_{\lambda} }   \dfrac{(qbr^{-1}x,
    qb/(axr^m))_{\mu} } {(x^{-1}, axr^{m-1})_{\mu} } \\  \cdot
    \dfrac{(r^m, qbr^{m-2})_{\lambda} }
    {(qbr^{-1}, r)_{\lambda} }  \dfrac{(r,
    br^{-1}t^{1-n})_{\mu}}{(qbr^{-2}, qt^{n-1})_{\mu}}
\cdot \prod_{i=1}^{n}\left\{ \dfrac{E(br^{-1}t^{2-2i} q^{2\mu_i})}
    {E(br^{-1}t^{2-2i})}  \left(qt^{2i-2} \right)^{\mu_i} \right\}\\
  \cdot\prod_{1\leq i< j\leq n}\left\{\dfrac{(qt^{j-i})_{\mu_i -\mu_j}
    } { (qt^{j-i-1})_{\mu_i - \mu_j} } \dfrac{
    (br^{-1}t^{3-i-j})_{\mu_i + \mu_j} } { (br^{-1}t^{2-i-j})_{\mu_i+
    \mu_j} } \right\} \\
\cdot W_{\mu} (q^{\lambda}t^{\delta(n)};q,p, t, br^{m-1}t^{2-2n},
    br^{-1}t^{1-n})
\end{multline}
where $xr^{\delta(m)} = (xr^{m-1},\ldots, xr, x)$ and $n$ is a
positive integer such that $n \geq \ell(\lambda)$.
\end{lemma}

\begin{proof}
The $m=1$ case of~(\ref{omega(xr^delta)}) reduces to the
definition~(\ref{eq:omega{lambda,mu}}). Assuming that the
identity~(\ref{omega(xr^delta)}) holds for all $k<m$, we expand the
\lhs\ of~(\ref{omega(xr^delta)}) using the recurrence
relation~(\ref{recurrence22}) and verify that it is summable by the
\ci~(\ref{CocycleIdentity}) giving the \rhs\ of~(\ref{omega(xr^delta)}).
Finally, an induction on $m$ gives the result to be proved.
\end{proof}

The next lemma will be crucial in proving the last two inductive steps
in Theorem~\ref{thm:bigInduction}.
\begin{lemma} Let $x,q, p, t,a,b\in \mathbb{C}$, then
  \label{lem:W-omega-rel}
  \begin{multline}
    \label{eq:W-omega-rel}
W_{\lambda/\mu}(x; q,p,t,a,b)\\ = \lim_{r\rightarrow t}\left\{\prod_{1\le
      i< j\le n}\left[\dfrac{(t^{j-i+1})_{\lambda_i
      -\lambda_j}(t^{j-i})_{\mu_i-\mu-j}(qbt^{2-i-j})_{\lambda_i
      +\lambda_j} (qbt^{-i-j})_{\mu_i+\mu_j}} {(t^{j-i})_{\lambda_i -
      \lambda_j} (t^{j-i+1})_{\mu_i-\mu_j} (qbt^{1-i-j})_{\lambda_i +
      \lambda_j} (qbt^{1-i-j})_{\mu_i +\mu_j}}\right]\right.\\
\left.\cdot\dfrac{(r)_\lambda(qbr^{-n})_\lambda (r^n)_\mu (qbr^{-2})_\mu}
      {(r^n)_\lambda (qbr^{-1})_\lambda (r)_\mu (qbr^{-n-1})_\mu}
      \omega_{\lambda/\mu}(x;r,a,b)\right\}
  \end{multline}
\end{lemma}
\begin{proof}
  We will prove Lemma \ref{lem:W-omega-rel} in three parts.  In the
  first part we reduce to proving the case when $\mu_n=0$.  If $\mu_n \ne 
0$,
  then we can express
\begin{equation}\label{W-mu-red}
W_{\lambda/\mu}(x;q,p, t,a,b)= f_1 \, W_{\hat{\lambda}/\hat{\mu}}
(xq^{-1}; q,p,t,aq^2,bq^2)
\end{equation}
and
\begin{equation}\label{omega-mu-red}
  \omega_{\lambda/\mu}(x;r, q, p, t; a,b)= f_2 \,
  \omega_{\hat{\lambda}/\hat{\mu}}
(xq^{-1}; r, q, p, t ;aq^2,bq^2)
\end{equation}
where $\hat{\lambda}=(\lambda_1 -1,\ldots,\lambda_n -1)$,
$\hat{\mu}=(\mu_1 -1,\ldots,\mu_n -1)$, and $f_1$,
$f_2$ are (complicated) factors.
A tedious computation shows that $f_1 = f_2$ and,
after iteration, we reduce to the  case $\mu_n = 0$.

Using the n-part Jackson sum (\ref{eq:W_Jackson}), we now will give an
explicit formula for $W_{\lambda/\mu}(t^{-n}; q,p,t,at^{2n}, bt^n)$ in
the case where $l(\lambda)= n$ and $l(\mu) < n$. For $s\in \mathbb{C}$
and $z =(z_1,\ldots, z_{n-1}) \in \mathbb{C}^{n-1}$, we consider
$$W_\lambda (s,z_1,\ldots, z_{n-1},1;q, p, t,a,b) = W_\lambda(s,z,1;q, p, 
t,a,b).$$
By identity (\ref{eq:ShiftedStable}), we have
\begin{equation}
  \label{eq:s-form-reduction}
  W_\lambda(s,z,1;q,p,t,a,b) =
  W_\lambda(st^{-1},t^{-1}z; q,p,t,at^{2n},bt^n).
\end{equation}
The property (\ref{def:item1}) for n-part partitions implies that
\begin{equation}
   \label{eq:s-symmetry}
  W_\lambda(st^{-1},t^{-1}z;q,p,at^{2n},bt^n)
  = W_\lambda(t^{-1}z, st^{-1};q,p,t,at^{2n},bt^n).
\end{equation}
Expanding the right-hand side of equation (\ref{eq:s-symmetry}), we
obtain
\begin{multline}\label{s-JacksonExpand}
  W_\lambda(r,z,1;q,p,t,a,b) \\
   = \sum_{\mu\subseteq
  \lambda}\dfrac{(s^{-1}t)_\lambda(ast^n)_\lambda}
  {(qbst^{n-1})_\lambda (qb/(as))_\lambda}\prod_{1\le i<j\le n}
  \dfrac{(t^{j-i+1})_{\lambda_i -\lambda_j}(qbt^{n+2-i-j})_{\lambda_i
  +\lambda_j}} {(t^{j-i})_{\lambda_i -\lambda_j} (qbt^{n+1-i-j})_{\lambda_i
  +\lambda_j}}\\ \cdot \dfrac{(bs)_\mu (qb/(at))_\mu} {(qt^{n-1})_\mu
  (ast^n)_\mu} \prod_{i=1}^n
  \dfrac{E(bst^{n+1-2i}q^{2\mu_i})(qt^{2i-2})^{\mu_i}}
  {E(bst^{n+1-2i})}\\ \cdot \prod_{1\le i<j\le n} \dfrac{(t^{j-i})_{\mu_i
  -\mu_j} (qt^{j-i})_{\mu_i -\mu_j}(qbst^{n-i-j})_{\mu_i +\mu_j}
  (bst^{n+2-i-j})_{\mu_i +\mu_j}} {(t^{j-i+1})_{\mu_i
  -\mu_j} (qt^{j-i+1})_{\mu_i -\mu_j}(qbst^{n+1-i-j})_{\mu_i +\mu_j}
  (bst^{n+1-i-j})_{\mu_i +\mu_j}}\\ \cdot W_\mu(q^\lambda
  t^{\delta(n)}; q,p,t,bt^{2-n},bs) W_\mu(s^{-1}z, 1; q,p,t,as^2,bs).
\end{multline}
We can express the last $W$ factor in expression
(\ref{s-JacksonExpand}) as
\begin{equation}
  \label{eq:st-reduction}
   W_\mu(s^{-1}z, 1; q,p,t,as^2,bs)= W_\mu((st)^{-1}z; q,p,t,a(st)^2,bst),
\end{equation}
which forces the partitions $\mu$ appearing in the sum
(\ref{s-JacksonExpand}) to have $l(\mu)<n$.

Using the the $W$-recurrence relation (\ref{eqWrecurrence}), we also
find
\begin{multline}
  \label{s-recurrence}
W_\lambda(s,z,p,1;q, p, t,a,b) \\ = \sum_{\mu\subseteq\lambda}
W_{\lambda/\mu}(st^{-n}; q,p,t,at^{2n},bt^n) W_\mu(z,1;q,p,t,a,b)\\ =
\sum_{\mu\subseteq\lambda} W_{\lambda/\mu}(st^{-n}; q,p,t,at^{2n},bt^n)
W_\mu(t^{-1}z; q,p,t,at^2,bt),
\end{multline}
where again  $l(\mu)<n$.
Taking the limits as $s\rightarrow 1$ of the right-hand sides of
identities (\ref{s-JacksonExpand}) and (\ref{s-recurrence}) and
comparing the coefficients of the (linearly independent) functions
$W_\mu(t^{-1}x; q,p,t,at^2,bt)$ for a fixed partition $\mu$, we obtain
the relation
\begin{multline}\label{eq:KeyLimit_s}
  W_{\lambda/\mu}(st^{-n}; q,p,t,at^{2n},bt^n)\\ = \lim_{s\rightarrow 1}
  \left\{\dfrac{(s^{-1}t)_\lambda(ast^n)_\lambda}
  {(qbst^{n-1})_\lambda (qb/(as))_\lambda}\prod_{1\le i<j\le n}
  \dfrac{(t^{j-i+1})_{\lambda_i -\lambda_j}(qbt^{n+2-i-j})_{\lambda_i
  +\lambda_j}} {(t^{j-i})_{\lambda_i -\lambda_j} (qbt^{n+1-i-j})_{\lambda_i
  +\lambda_j}}\right.\\ \left.\cdot \dfrac{(bs)_\mu (qb/(at))_\mu}
{(qt^{n-1})_\mu (ast^n)_\mu} \prod_{i=1}^n
  \dfrac{E(bst^{n+1-2i}q^{2\mu_i})(qt^{2i-2})^{\mu_i}}
  {E(bst^{n+1-2i})}\right.\\ \left.\cdot \prod_{1\le i<j\le n}
\dfrac{(t^{j-i})_{\mu_i
  -\mu_j} (qt^{j-i})_{\mu_i -\mu_j}(qbst^{n-i-j})_{\mu_i +\mu_j}
  (bst^{n+2-i-j})_{\mu_i +\mu_j}} {(t^{j-i+1})_{\mu_i
  -\mu_j} (qt^{j-i+1})_{\mu_i -\mu_j}(qbst^{n+1-i-j})_{\mu_i +\mu_j}
  (bst^{n+1-i-j})_{\mu_i +\mu_j}}\right.\\ \left.\cdot W_\mu(q^\lambda
  t^{\delta(n)}; q,p,t,bt^{2-n},bs) \right\}.
\end{multline}
Substitute $r=s^{-1}t$ into relation (\ref{eq:KeyLimit_s}) and
take the limit as $r\rightarrow t$ in place of $s\rightarrow 1$. Using
the definition of the $\omega_{\lambda/\mu}$ function given by
  identity (\ref{eq:omega{lambda,mu}}), then one observes that we can
replace the  right-hand side of (\ref{eq:KeyLimit_s}) to find
  \begin{multline}\label{SpecialLimitOmega}
W_{\lambda/\mu}(st^{-n}; q,p,t,at^{2n},bt^n)\\
=   \lim_{r\rightarrow t}\left\{\prod_{1\le
      i< j\le n}\left[\dfrac{(t^{j-i+1})_{\lambda_i
      -\lambda_j}(t^{j-i})_{\mu_i-\mu-j}(qbt^nt^{2-i-j})_{\lambda_i
      +\lambda_j} (qbt^nt^{-i-j})_{\mu_i+\mu_j}} {(t^{j-i})_{\lambda_i -
      \lambda_j} (t^{j-i+1})_{\mu_i-\mu_j} (qbt^nt^{1-i-j})_{\lambda_i +
      \lambda_j} (qbt^nt^{1-i-j})_{\mu_i +\mu_j}}\right]\right.\\
\left.\cdot\dfrac{(r)_\lambda(qbt^nr^{-n})_\lambda (r^n)_\mu 
(qbt^nr^{-2})_\mu}
      {(r^n)_\lambda (qbt^nr^{-1})_\lambda (r)_\mu (qbt^nr^{-n-1})_\mu}
      \omega_{\lambda/\mu}(t^{-n}; r; at^{2n},bt^n)\right\}
  \end{multline}

Identity (\ref{SpecialLimitOmega})is simply  the $x=t^{-n}$ (and
$l(\mu)<n$) case of Lemma \ref{lem:W-omega-rel}.  We complete the
proof of Lemma \ref{lem:W-omega-rel} by observing that if Lemma
\ref{lem:W-omega-rel} is true for this special value of $x$, then it must
be true for all $x$ in the domains of $W_{\lambda/\mu}$ and
$\omega_{\lambda/\mu}$. This is because the factors involving the
variable $x$ in $W_{\lambda/\mu}(x; q,p,t,at^{2n},bt^n)$ by its
definition (\ref{definitionSkewW}) and in
$\omega_{\lambda/\mu}(t^{-n};r; at^{2n},bt^n)$ by  its definition
(\ref{eq:omega{lambda,mu}}) are identical in the limit $r\rightarrow
1$. This means that the two sides of (\ref{eq:W-omega-rel}) agree up
to a factor independent of $x$, which must be one. This completes the
proof of Lemma \ref{lem:W-omega-rel}.
\end{proof}

Applying the $W$-recurrence relation (\ref{eqWrecurrence}) together
with the corresponding $\omega$-recurrence relation
\ref{recurrence22}, we obtain an immediate corollary of Lemma
\ref{lem:W-omega-rel}:
\begin{cor} Let $k$ be a positive integer and $z=(z_1,\ldots, z_k)\in
  \mathbb{C}^k$. Let $q, p, t,a,b\in \mathbb{C}$, then
  \label{cor:W-omega-rel}
  \begin{multline}
    \label{coreq:W-omega-rel}
W_{\lambda/\mu}(z; q,p,t,a,b)\\ = \lim_{r\rightarrow t}\left\{\prod_{1\le
      i< j\le n}\left[\dfrac{(t^{j-i+1})_{\lambda_i
      -\lambda_j}(t^{j-i})_{\mu_i-\mu-j}(qbt^{k+1-i-j})_{\lambda_i
      +\lambda_j} (qbt^{-i-j})_{\mu_i+\mu_j}} {(t^{j-i})_{\lambda_i -
      \lambda_j} (t^{j-i+1})_{\mu_i-\mu_j} (qbt^{k-i-j})_{\lambda_i +
      \lambda_j} (qbt^{1-i-j})_{\mu_i +\mu_j}}\right]\right.\\
\left.\cdot\dfrac{(r)_\lambda(qbr^{-n+ k-1})_\lambda (r^n)_\mu
      (qbr^{-2})_\mu}
      {(r^n)_\lambda (qbr^{k-2})_\lambda (r)_\mu (qbr^{-n-1})_\mu}
      \omega_{\lambda/\mu}(z; r ;a,b)\right\}
  \end{multline}
\end{cor}

Setting $z= (x,y)$ or $z = (y,x)$ for $x,y \in \mathbb{C}$ in
identity (\ref{coreq:W-omega-rel}), and using
the symmetry of $\omega_{\lambda/\mu}(x,y;r,a,b)$ in the variables
$x$ and $y$,
it follows that $W_{\lambda/\mu}(x,y; q, p, t,a,b)$ is also symmetric in the
variables $x$ and $y$.  This proves step seven of  Theorem
\ref{thm:bigInduction}.

To prove step eight, we will consider the function
$W_{\lambda}(z,t^{\delta(k)};q, p, t,a,b)$,
for $z\in \mathbb{C}^{n+1}$, $l(\lambda)\le n$ and positive integer
k. It follows from identities (\ref{eq:ShiftedStable}),
(\ref{eqWrecurrence})
and step seven, that
\begin{multline}\label{key_expandJack}
W_{\lambda}(zt^{-k};q, p, t,at^{2k},bt^k) \\=
W_{\lambda}(z,t^{\delta(k)};q, p, t,a,b) =
W_{\lambda}(t^{\delta(k)},z;q, p, t,a,b)\\ = \sum_{\mu\subseteq \lambda}
W_{\lambda/\mu}(t^{-n-1} t^{\delta(k)};q, p, t,at^{2n+2}, bt^{n+1}) \cdot
W_{\mu}(z; q, p, t,a,b).
\end{multline}
Substituting identities (\ref{coreq:W-omega-rel}) and
(\ref{omega(xr^delta)}) into (\ref{key_expandJack}) and simplifying,
we obtain the ``extended'' $W$-Jackson sum (\ref{eq:W_Jackson}) in the
case of a (n=1)-tuple variable $x$ and $l(\lambda)\le n$.  This
completes the proof of step eight and also Theorem \ref{thm:bigInduction}.
\end{proof}

\section{$\omega$--\Js}
\label{section4}
In this section we prove an extension of the
$W$--\Js~(\ref{eq:W_Jackson}) and obtain further properties of the
\Jc. First note that  
the symmetry (Lemma~\ref{symmetry}) of the infinite matrix \W{x_1,x_2} in 
$x_1,x_2\in\mathbb{C}$ implies that
\begin{equation}
  \omega_\lambda(1,x_1; r; a,b) = \W{1,x_1}_{\lambda0} =
  \W{x_1,1}_{\lambda0}
\end{equation}
Entry level writing of this relation gives
\begin{multline}
\label{JacksonOnedim}
\omega_\lambda(1,x_1; r; a,b) = \sum_\mu \omega_{\lambda/\mu}(r^{-1}; r;
  ar^{2}, b r ) \, \omega_{\mu}(x_1; r; a, b) \\ = \sum_\mu
  \omega_{\lambda/\mu}(r^{-1}x_1; r; ar^{2}, b r ) \, \omega_{\mu}(1; r; a,
  b) = \omega_{\lambda}(r^{-1}x_1; r; ar^2, br)
\end{multline}
where the last equality follows from the obvious identity
\begin{equation}
\label{vanishes}
  \omega_{\mu}(1; r; a, b) = \delta_{\mu 0}
\end{equation}

In other words, setting $x_2=1$ in $\omega_\lambda(x_2,x_1;a,b)$ will
have the effect of shifting the parameters.  The
identity~(\ref{JacksonOnedim}) simplifies to give the summation
formula
\begin{equation}
\label{JsinTermsOfw}
  \omega_\lambda(r^{-1}x_1; r; ar^2,bu) =\sum_\mu
  \omega_{\lambda/\mu}(r^{-1}; r; ar^2,br) \, \omega_\mu(x_1; r; a,b)
\end{equation}
which is exactly the one variable version of the
$W$--\Js~(\ref{eq:W_Jackson}) 
\begin{multline}
\label{simplifiedJs2}
  \dfrac{(sx^{-1}, as x)_\lambda}{(qb x, q b /ax )_\lambda}
 =\sum_{\mu \subseteq \lambda} \dfrac{(s, a
   s)_\lambda} {(qb, qb /a)_\lambda} \dfrac{ (bt^{1-n}, qb/a s)_\mu }{
   (qt^{n-1}, as)_\mu} \\ \cdot\prod_{i=1}^n \left\{\dfrac{E(b
   t^{2-2i}q^{2\mu_i})} {E(b t^{2-2i})} (qt^{2i-2})^{\mu_i}\right\}
 \prod_{1\leq i<j \leq 
 n} \left\{\dfrac{ (qt^{j-i})_{\mu_i-\mu_j} (bt^{3-i-j})_{\mu_i+\mu_j}}
 {(qt^{j-i-1})_{\mu_i-\mu_j} (b t^{2-i-j})_{\mu_i+\mu_j} }
 \right\}  \\ \cdot W_\mu(q^\lambda t^{\delta(n)}; q,p, t, b rt^{2-2n},
 bt^{1-n}) \dfrac{(x^{-1}, ax)_\mu} {(qbx, q b /ax)_\mu} 
\end{multline}
after a simple
reparametrization. We expect to have a similar
formulation of the \Js\ in the multivariable case.

Before we proceed to write a multivariable version of the \Js\ in
terms of $\omega_{\lambda/\mu}$ function, we verify the generalization
of the one variable result~(\ref{vanishes}) in
\begin{lemma}
\label{delta{lambda,0}}
With the notation as above and for any partition $\lambda$, we have
\begin{equation}
   \omega_{\lambda}(r^{\delta(m)}; r; a, b) = \delta_{\lambda 0}
\end{equation}
where $r^{\delta(m)} = (r^{m-1},\ldots, r,1)$ for any
$m\in\mathbb{Z}_{\geq}$.
\end{lemma}

\begin{proof}
One variable (i.e., $m=1$) case $\omega_{\lambda}(1;r; a, b) =
  \delta_{\lambda 0}$ is obvious. Assume that the identity holds for
  all $k<m$, that 
  is $\omega_{\lambda}(r^{\delta(k)}; r; a, b) = \delta_{\lambda 0}$
  for such $k$. We'd like to show
  that it also holds true for $m$.  But,
\begin{multline}
\omega_{\lambda}(r^{\delta(m)}; r; a, b) = \W{r^{\delta(m)} }_{\lambda0}
  \\ = \sum_{\mu} \W{1}_{\lambda\mu} \, \W{r^{\delta(m-1)} }_{\mu0} =
  \W{1; r; ar^{2m-2}, br^{m-1}}_{\lambda0}
\end{multline}
by the induction hypothesis. Therefore,
\begin{equation}
\omega_{\lambda}(r^{\delta(m)}; r; a, b) = \omega_{\lambda}(1; r; ar^{2m-2},
br^{m-1}) =\delta_{\lambda 0}
\end{equation}
due to the obvious one variable result.
\end{proof}

We now generalize the argument used above in one variable case to
write a multivariable version of the \Js\ in
\begin{lemma}
\label{multiJs}
Let $k,m\in\mathbb{Z}_{\geq}$, $r\in\mathbb{C}$ and $z=(x_1,\ldots,
x_k)\in\mathbb{C}^k$, and let $r^{-m}z$ denote
$r^{-m}z=(r^{-m}x_1,\ldots, r^{-m}x_k)$ as before.  Then, we have
\begin{equation}
\label{eqn:multiJs11}
  \omega_{\lambda}(r^{-m}z; r; ar^{2m}, br^m) =
  \omega_\lambda(r^{\delta(m)},z; r; a,b)
\end{equation}
where $\lambda$ is any partition.
\end{lemma}

\begin{proof}
This could be seen from the symmetry property of $\omega_\lambda$
which is shown to be equivalent to the fact that $\mathfrak{W}(x_i)$
commute. For any $m\in\mathbb{Z}_{>}$, we have
\begin{equation}
  \omega_\lambda(r^{\delta(m)},z; r; a,b) = \W{r^{\delta(m)},z}_{\lambda0}
  = \W{z,r^{\delta(m)}}_{\lambda0}
\end{equation}
which may be expanded as
\begin{multline}
\label{eqn:multiJs}
\omega_\lambda(r^{\delta(m)},z; r; a,b) = \sum_\mu
  \omega_{\lambda/\mu}(r^{-k}r^{\delta(m)}; r; ar^{2k}, b r^k) \,
  \omega_{\mu}(z; r; a, b) \\ = \sum_\mu \omega_{\lambda/\mu}(r^{-m}z;
  r; ar^{2m}, b r^m ) \, \omega_{\mu}(r^{\delta(m)}; r; a, b) =
  \omega_{\lambda}(r^{-m}z; r; ar^{2m}, br^m)
\end{multline}
where the last equality follows from Lemma~\ref{delta{lambda,0}}.
\end{proof}

Using the identity~(\ref{eqn:multiJs}) and~(\ref{omega(xr^delta)}) we
can now write the multivariable \Js\ explicitly. We have
\begin{theorem}[Multivariable $\omega$--\Js]
\label{thm:multiJsExplicit}
For a partition $\lambda$,
and complex parameters $q,p,t,r,s,a,b\in\mathbb{C}$, we have
\begin{multline}
\label{10phi9sum}
\dfrac{ (qb/a, qb)_{\lambda}}{ (as, s)_{\lambda} } \dfrac{
    \omega_{\lambda}(s^{-1}z; r; ar^{1-k}s^{2}, br^{1-k}s)} {
    \omega_{\lambda}(r^{-k}; r; qbsr^{k-1}, bs) } \\
= \sum_\mu \prod_{i=1}^{n}\left\{ \dfrac{E(b t^{2-2i} q^{2\mu_i})}
    {E(b t^{2-2i})}  \left(qt^{2i-2} \right)^{\mu_i} \right\}
  \hspace{-5pt} \prod_{1\leq i< j \leq n} \left\{ \dfrac{
    (qt^{j-i})_{\mu_i - \mu_j} } { (qt^{j-i-1})_{\mu_i - \mu_j} } \dfrac{
    (b t^{3-i-j})_{\mu_i + \mu_j} } { (b t^{2-i-j})_{\mu_i +
    \mu_j} } \right\} \\   \cdot
\dfrac{ (b t^{1-n}, qb/as)_{\mu}}{ (qt^{n-1}, as )_{\mu}} \,
W_{\mu} ( q^{\lambda}t^{\delta(n)};q, p, t, bs t^{2-2n},
b t^{1-n})  \, \dfrac{ \omega_{\mu}(z; r; ar^{1-k}, br^{1-k}) } {
    \omega_{\mu}(r^{-k}; r; qbr^{k-1}, b) }
\end{multline}
where $z=(x_1,\ldots,x_k)\in\mathbb{C}^k$ for some
$k\in\mathbb{Z}_{>}$  and $n$ is a positive integer such that $n \geq
\ell(\lambda)$.
\end{theorem}

\begin{proof}
It follows by induction from the Lemma~(\ref{multiJs})
and the key identity~(\ref{omega(xr^delta)}) that
\begin{multline}
\label{multiJsExplicit}
\dfrac{(qbr^{m-1}, qb/a, qbr^{k-1}, r)_{\lambda}
    } {(r^k, ar^{k+m-1}, r^m, qbr^{k+m-2})_{\lambda} }
\, \omega_{\lambda}(r^{-m}z; r;ar^{2m}, br^m)  \\
= \sum_\mu  \prod_{i=1}^{n}\left\{ \dfrac{E(br^{k-1}t^{2-2i} q^{2\mu_i})}
    {E(br^{k-1}t^{2-2i})}  \left(qt^{2i-2} \right)^{\mu_i} \right\}
  \dfrac{(r, br^{k-1}t^{1-n})_{\mu}}{(qbr^{k-2},
qt^{n-1})_{\mu}} \\ \cdot\dfrac{(qbr^{-1}, qb/ar^{m})_{\mu} } {(r^k,
ar^{k+m-1} )_{\mu} } \prod_{1\leq i< j
    \leq n} \left\{ \dfrac{ (qt^{j-i})_{\mu_i - \mu_j}
    } { (qt^{j-i-1})_{\mu_i - \mu_j} } \dfrac{
    (br^{k-1}t^{3-i-j})_{\mu_i + \mu_j} } { (br^{k-1}t^{2-i-j})_{\mu_i +
    \mu_j} } \right\} \\  \cdot
W_{\mu} ( q^{\lambda}t^{\delta(n)};q, p, t, br^{k+m-1}t^{2-2n},
br^{k-1}t^{1-n})  \, \omega_{\mu}(z;r; a, b)
\end{multline}
where $z=(x_1,\ldots,x_k)\in\mathbb{C}^k$ and $n$ is a
positive integer such that $n \geq \ell(\lambda)$.

Both sides of this identity are
rational functions in the single parameter $s$. We rewrite this
identity
with some elementary algebra as an
identity between two polynomials that are equal at infinitely many
values $s=r^m$ for any positive integer $m$
by~(\ref{multiJsExplicit}).
\end{proof}

\begin{remark}
The identity~(\ref{10phi9sum}) looks like a \tns\ summation identity,
for it has an extra free parameter $r\in\mathbb{C}$. However, it
reduces to $W$--\Js~(\ref{eq:W_Jackson}) if $z=x_1\in\mathbb{C}$ or in
the case when we send $r\rightarrow t$ and set
$k=n$. Therefore, the identity~(\ref{10phi9sum}) extends $W$--\Js\ 
to the case where $r\in\mathbb{C}$ and $k\in\mathbb{Z}_{>}$ are
arbitrary. 
\end{remark}

\begin{remark}
Besides extending it, the preceeding argument also provides an
alternative proof to $W$--\Js, for the proof of the
$\omega$--\Js~(\ref{10phi9sum}) depends only on
the \ci~(\ref{CocycleIdentity}) which can be independently proved by
other methods (see the bulk difference equation in~\cite{Rains1}, for
example). 
\end{remark}

We conclude this section by listing some important properties of the
$\omega_{\lambda/\mu}$ functions.
\begin{enumerate}
\item $\omega_{\lambda}$ satisfy the vanishing property
\begin{equation}
\label{omegaFundamentalVanishing}
  \omega_{\lambda}(q^\pi r^{\delta(n)};r, q, p,t; a,b)=0
\end{equation}
where $\pi=(\pi_1,\ldots,\pi_n)$ and $\lambda
=(\lambda_1,\ldots,\lambda_n)$ are $n$-part partitions such that
$\lambda \nsubseteq \pi_1^n$. The proof is a slight generalization of
that of Lemma~(\ref{delta{lambda,0}}).
\item Setting $z=r^{\delta(k)}$ in~(\ref{10phi9sum}) gives
\begin{equation}
  \omega_{\lambda}(xr^{\delta(k)}; r; a, b) = \dfrac{ (x^{-1},
    ar^{k-1}x )_{\lambda} }  { (qbr^{k-1}x, qb/ax)_{\lambda}} 
   \dfrac{(r^k, qbr^{k-2})_\lambda }{(qbr^{-1}, r)_\lambda } 
\end{equation}
\item If we set $s=r^{-1}$ in the
  \ci~(\ref{CocycleIdentity}) we get the following inversion relation
\begin{equation}
\delta_{\nu\tau} = \sum_{\tau\subseteq
  \lambda\subseteq \nu} \omega_{\nu/\lambda}(r;r^{-1}; a, b)
\,  \omega_{\lambda/\tau}(r^{-1};r; ar^2,br)
\end{equation}
where $\delta_{\nu\tau} = 1$ when $\nu=\tau$, and $\delta_{\nu\tau} =
0$ otherwise. The proof is a special case of the proof of
Lemma~(\ref{cor:CocycleIdentity}).
\item We have the following elliptic transformation identities: For
  $x\in\mathbb{C}$
\begin{multline}
\label{ellipticInb}
\omega_{\lambda/\mu}(x; r; a, pb) 
= (q^2 b^2 /a)^{|\lambda|-|\mu|} t^{-2n(\lambda)+2n(\mu)} 
  q^{2n(\lambda') -2n(\mu')} \omega_{\lambda/\mu}(x; r; a, b) 
\end{multline}
\begin{multline}
\label{ellipticIna}
\omega_{\lambda/\mu}(x; r; pa, b) \\ =
(q b)^{-|\lambda|+|\mu|} p^{|\lambda|-|\mu|} r^{-|\mu|}
t^{2n(\lambda)- 2n(\mu)} q^{-2n(\lambda')+2n(\mu')} 
\omega_{\lambda/\mu}(x; r; a, b) 
\end{multline}
\begin{equation}
\label{ellipticInr}
\omega_{\lambda/\mu}(x; pr; a, b) 
= (ar^{-2})^{-|\mu|} p^{2|\mu|}
  t^{2n(\mu)} q^{-2n(\mu')} \omega_{\lambda/\mu}(x; r; a, b) 
\end{equation}
\begin{equation}
\label{ellipticInx}
\omega_{\lambda/\mu}(px; r; a, b) 
= \omega_{\lambda/\mu}(x; r; a, b) 
\end{equation}
Proofs use direct computation and properties like  
\begin{equation}
  (p^{-1}x)_{\lambda} = (-1)^{|\lambda|} p^{-|\lambda|} x^{|\lambda|}
  t^{-n(\lambda)} q^{n(\lambda')} (x)_\lambda 
\end{equation}
and
\begin{equation}
  (px)_{\lambda} = (-1)^{|\lambda|} x^{-|\lambda|} t^{n(\lambda)}
  q^{-n(\lambda')} (x)_\lambda.
\end{equation}
Note that~(\ref{ellipticInb}) and~(\ref{ellipticInx}) of
these identities can be generalized to the multivariable case
$x\in\mathbb{C}^n$ using
the recurrence formula~(\ref{recurrence22}) for $\omega_{\lambda/\mu}$. 
\item For $\mu=\lambda$ we get
\begin{multline}
\omega_{\lambda/\lambda}(x; r; a, b)
= \dfrac{(qbr^{-1}x, qb/arx)_\lambda }{(qbx, qb/ax)_\lambda
} \dfrac{(br^{-1}t^{1-n}, qbr^{-1} )_{\lambda}} {(bt^{1-n},
  qbr^{-2})_{\lambda}} \\ \cdot r^{|\lambda|}  
 \prod_{i=1}^{n}\left\{ \dfrac{(bt^{2-2i})_{2\lambda_i} }
    {(br^{-1}t^{2-2i})_{2\lambda_i} } \right\} 
\prod_{1\leq i< j \leq n} \left\{
   \dfrac{(bt^{2-i-j})_{\lambda_i+\lambda_j}} 
{(bt^{3-i-j})_{\lambda_i+\lambda_j}} \dfrac{
    (br^{-1}t^{3-i-j})_{\lambda_i + \lambda_j} } {
    (br^{-1}t^{2-i-j})_{\lambda_i + \lambda_j} } \right\}
\end{multline}
The proof follows from Proposition~\ref{specialevalW}.
\end{enumerate}
Other properties including a reversal, an inversion and a duality
identity for $\omega_{\lambda/\mu}$ will be given in a later
publication.  

\bibliographystyle{amsalpha}

\begin{thebibliography}{ABC}
\bibitem[AAB]{AgarwalA1}     A. K.   Agarwal, G. E. Andrews and
  D. M. Bressoud,      \emph{The Bailey lattice,}      J. Indian
  Math. Soc.            \textbf{51}     (1987), 57--73. 
\bibitem[AAR]{AndAskRoy}  G. E. Andrews, R. Askey and R. Roy
    \emph{Special Functions, Encyclopedia of mathematics and its
  applications,} Cambridge  University Press, Cambridge,
  (1999). 
\bibitem[ASW]{AndrewsS1}     G. E.   Andrews, A. Schilling and
  S. O. Warnaar,        \emph{An $A_2$ Bailey lemma and
  Rogers--Ramanujan--type identities,}       J. Amer. Math. Soc.
  \textbf{12}     (1999), 677--702. 
\bibitem[An]{Andrews}  G. E. Andrews \emph{q-Series:  Their Development
    and Application in Analysis, Number Theory, Combinatorics, Physics
    and Computer Algebra.}  NSF CBMS Regional Conference Series,
    \textbf{vol. 66}. (1986), Washington, DC: CBMS
\bibitem[B1]{Bailey1}       W. N.   Bailey, \emph{Identities of the
  Rogers--Ramanujan type,} Proc. London Math. Soc. (2)     \textbf{50}
  (1949), 1--10. 
\bibitem[B2]{Bailey2}    \bysame, \emph{Some identities in
  combinatory analysis,} Proc. London Math. Soc. (2)     \textbf{49}
  (1947), 421--435. 
\bibitem[BL]{BiedenharnL1}  L. C.   Biedenharn and J. D. Louck,
  \emph{A new class of symmetric polynomials defined in terms of
  tableaux,}       Adv. in Appl. Math.             \textbf{10}
  (1989), 396--438.
\bibitem[Br]{Bressoud1}  D. M.   Bressoud,       
\emph{The Bailey
  lattice, an introduction,} in Ramanujan revisited, (G. E. Andrews
  et al. eds.), Academic Press, NewYork,              \textbf{}
  (1988), 57--67. 
\bibitem[C1]{Coskun} H. Coskun,
\emph{A $BC_n$ Bailey Lemma and generalizations of \RRis,} August
  2003, Ph.D. thesis.
\bibitem[C2]{Coskun2} \bysame , % H. Coskun,
\emph{Elliptic $BC_n$ Bailey Lemmas and applications,} preprint.
\bibitem[CG]{Coskun1}  H. Coskun and R. Gustafson,
\emph{$BC_n$ Bailey lemmas and generalizations of Rogers--Ramanujan
  identities,}  preprint.
\bibitem[FT]{FrenkelT1}     L. B.   Frenkel and V. G. Turaev,
  \emph{Elliptic solutions of the Yang--Baxter equation and modular
  hypergeometric functions,}    The Arnold--Gelfand Mathematical
  Seminars. Birkhauser, Boston, MA   \textbf{}
  (1997), 171--204.
\bibitem[GIS]{GarretI1}       K.      Garret, M. E. H. Ismail and
  D. Stanton, \emph{Variants of the \RRis,}   Adv. in App. Math. 
  \textbf{23} (1999), 274--299.    
\bibitem[GR]{GasperR1}      G.      Gasper and M. Rahman,   \emph{Basic
  hypergeometric series, Encyclopedia of mathematics and its
  applications,} Vol 35 Cambridge  University Press, Cambridge,
  (1990).
\bibitem[GH]{GouldenH1}     I. P.   Goulden and A. M. Hamel,
  \emph{Shift operators and factorial symmetric functions,}
  J. Combin. Theory Ser. A (1) \textbf{69}     (1995), 51--60.
\bibitem[Kn]{Knop1}   F. Knop, \emph{Symmetric and nonsymmetric quantum
    Capelli polynomials}, Comment. Math. Helv. \textbf{72} (1997),
  84--100. 
\bibitem[KS]{KnopSahi}         F.   Knop and S. Sahi,
  \emph{Difference equations and symmetric polynomials defined by
  their zeros}, Internat. Math. Res. Notices  \textbf{10} (1996), 473--486
\bibitem[LM]{LillyM3}       G. M.   Lilly and S. C. Milne,  \emph{The
  $A_\ell$ and $C_\ell$ Bailey transform and lemma,}
  Bull. Amer. Math. Soc. (N. S.)          \textbf{26}     (1992),
  258--263.
\bibitem[M1]{Macdonald1}    I. G.   Macdonald,
 \emph{Schur functions:
   theme and variations, Seminaire Lotharingien de combinatoire,}
   Publ. IRMA, Univ. Louis Pasteur, 498/S--27, Strasbourg,
   \textbf{}       (1992), 5--39.
\bibitem[M2]{Macdonald5}    I. G.   Macdonald,
\emph{Symmetric
  functions and Hall Polynomials, 2nd ed.,}  Oxford
  University Press, New York             \textbf{} 1995.
\bibitem[O1]{Okounkov1}     A.      Okounkov,       \emph{On Newton
  interpolation of symmetric functions: A characterization of
  interpolation Macdonald Polynomials,}       Adv. in Appl. Math.
  \textbf{20}     (1998), 395--428.
\bibitem[O2]{Okounkov2}     A.      Okounkov,    \emph{BC-Type
    Interpolation Macdonald Polynomials and Binomial Formulas for
    Koornwinder polynomials,}      Transform. Groups, \textbf{3(2)}
    (1998), 181--207 
\bibitem[OO]{OkOl} A. Okounkov and G. Olshanski, Shifted Schur
 functions, Algebra i Analiz \textbf{9} (1997), No. 2 (Russian)
\bibitem[R1]{Rains1} E. Rains, \emph{$BC_n$--symmetric abelian
    functions,} math.CO$/$0402113
\bibitem[HR]{Rosengren1}       H.   Rosengren, \emph{Elliptic
  hypergeometric series on root systems,}  Adv. Math. \textbf{181}
  (2004), 417-447. 
\bibitem[Sa]{Sahi1}  S. Sahi, Interpolation, integrality and a
  generalization of 
Macdonald's polynomials, Int. Math. Res. Not. \textbf{10} 1996, 457-471
\bibitem[SZ]{SpiridonovZ1}  V.      Spiridonov and A. Zhedonov,
  \emph{Classical biorthogonal rational functions on elliptic grids,}
  C. R. Math. Rep. Acad. Sci. Canada ,            \textbf{22}
  (2000), 70--76. 
\bibitem[Sp]{Spiridonov1}  V.      Spiridonov,
  \emph{An elliptic incarnation of the Bailey chain,}
  Internat. Math. Res. Notices 37 (2002), 1945-1977. 
\bibitem[W1]{Warnaar1}      S. O.   Warnaar,        \emph{50 years of
  Bailey's lemma,}      Kerber Festschrift. %,             \textbf{}      
\bibitem[W2]{Warnaar2}      S. O.   Warnaar,        \emph{Summation and
  transformation formulas for elliptic hypergeometric series,}
  Constr. Approx.         \textbf{18}     (2002), 479--502.
\end{thebibliography}

\end{document}